\documentclass{book}
\usepackage{amsmath}
\usepackage{amssymb}
\usepackage{multirow}
\usepackage{amsthm}
\usepackage{natbib}
\usepackage[bookmarks=false]{hyperref}
\usepackage{enumerate}

\input xy\huge
\xyoption{all}

\usepackage[pdftex]{graphicx}

\newtheorem{theorem}{Theorem}[chapter]

\newtheorem{lemma}[theorem]{Lemma}

\newtheorem{corollary}[theorem]{Corollary}

\theoremstyle{definition}
\newtheorem{remark}[theorem]{Remark}
\newtheorem*{remark*}{Remark}
\newtheorem{definition}[theorem]{Definition}
\newtheorem*{definition*}{Definition}

\newcommand{\rr}{\mathbb{R}}
\newcommand{\cc}{\mathbb{C}}

\newcommand{\cl}{\mathcal{C}}

\newcommand{\mm}{\mathcal{M}}
\newcommand{\kk}{\mathbb{K}}

\newcommand{\hh}{\mathbb{H}}

\newcommand{\belt}{\mbox{Belt}}

\newcommand{\End}{\mbox{End}}
\newcommand{\diff}{\mbox{Diff}_+}

\newcommand{\jj}{\mathcal{J}}

\begin{document}
%\maketitle
\begin{titlepage}
\begin{center}
\begin{Large}
Hebrew University of Jerusalem\\
\vspace{5mm}
Faculty of Science\\
\vspace{5mm}
Einstein Institute of Mathematics
\end{Large}
\end{center}

\vspace{20mm}
\begin{center}
\begin{huge}
An Invitation to \\
\vspace{4mm}

Topological Conformal Field Theories
\end{huge}
\end{center}

\vspace{30mm}
\begin{center}

\begin{large}
Submitted by
\end{large}

\vspace{5mm}
\begin{Large}
Amitai Zernik
\end{Large}

\vspace{5mm}
\begin{large}
guided by
\end{large}

\vspace{5mm}
\begin{Large}
Prof. David Kazhdan
\end{Large}
\end{center}

\begin{large}
\noindent
as part of the requirements for the degree of M. Sc. in Mathematics.
\end{large}

\vspace{40mm}

\begin{center}
\begin{Large}
September 2010
\end{Large}
\end{center}

\end{titlepage}

\chapter*{Preface}

This work is submitted as an MSc thesis to the Hebrew University of Jerusalem.

I have been fortunate to have David Kazhdan as my advisor. I've learned some beautiful mathematics from him. Better still, he has taught me how to think - mostly by serving an excellent example.

I am also grateful to Jake Solomon and Genadi Levin for a number of fruitful discussions.

\tableofcontents
\newpage
\chapter{Introduction}
In this work we present what seems to be a novel approach to conformal field theories in the spirit of \citep{segal}. Compared with the classical approach to the subject, our definition relies on only very few results from the theory of quasiconformal mappings (most necessary results are summarized in chapter \ref{sec:QCmappings}). This is an obvious advantage to those wanting to learn or teach the subject. In fact, this work grew out of our own attempt to study \citep{costello}, and we will mostly adhere to the notations and terminology therein.

Conformal field theories should be thought of as representations of Segal's category, a category enriched over topological spaces. It is a kind of cobordism category: the objects of this category are intervals and circles, and the morphisms are represented by Riemann surfaces with boundary with the in and out objects embedded into the boundary of the surface.

Our first task will be to define Segal's category. As we will see, defining the set of morphisms will be quite easy. It is the definition of the topology of these sets (``when are two Riemann surfaces close to each other?'') and the composition of morphisms (``how can two Riemann surfaces be sewed?'') which will require quasiconformal mappings.

Quasiconformal mappings can be thought of as mappings which have the right amount of flexibility: they are flexible enough so that we may freely deform a surface, thus obtaining topologically equivalent (but conformally distinct) surfaces, and also flexible enough so that we may align the boundaries of two surfaces, allowing them to be sewed. Yet these maps retain much of the appealing regularity of conformal mappings, allowing us to economically keep track of all this stretching and deforming by means of Beltrami differentials (which are locally just measurable complex-valued functions).

Our work is structured as follows. In chapter \ref{chapt2} we define the objects and the morphism sets of Segal's category. In chapter \ref{chapt3} we make a didactic attempt to topologize the morphism sets and define composition maps using smooth mappings. This attempt will present us with ample evidence to the necessity of a more flexible kind of mapping, thus motivating the introduction of quasiconformal mappings in the next chapter. Chapter \ref{chapt4} surveys, with accurate references but no proofs, some results from the theory of QC mappings and also contains proofs of corollaries which facilitate the application of these results.

Chapters \ref{chapt5} and \ref{chapt6} contain the heart of our work. In chapter \ref{chapt5} we define quasiconformal surfaces and the bundle of Beltrami differentials over it. We can then represent morphisms as measurable sections of this bundle. In chapter \ref{chapt6} we use this improved representation to define the topology of the morphism spaces and the composition maps, and prove that they are continuous and associative.

Chapter \ref{chapt7} and \ref{chapt8} fill in some details in the definition of conformal field theories and then define \emph{topological} conformal field theories. The latter are also representations of a category; this category is obtained from Segal's category by replacing the morphisms spaces by chain complexes and adjusting the composition maps accordingly.

Appendix \ref{AppA} contains a proof of a classification theorem for parameterized surfaces, based on the well-known classification of surfaces with boundaries. In appendix \ref{AppB} we prove that the image of the boundary after the sewing is smooth.

\chapter{The Objects and Morphisms of Segal's Category}
\label{chapt2}

We will denote Segal's category by $\mathcal{M}_\Lambda$. It is a symmetric monoidal category enriched over the category of topological spaces (cf. chapter \ref{sec:cft-and-tcft}).

Here $\Lambda$ represents a certain set, called ``the set of D-branes''. To ease the exposition and notation we assume this set consists of a single element. The simple modifications needed to accommodate $|\Lambda| > 1$ will be explained in chapter \ref{sec:finishing}.

\vspace{2mm}
An object of $M_\Lambda$ is a pair of non-negative integers, $(C,O)$. One can think of this as short hand for $(S^1)^C \coprod [0,1]^O$ with some fixed orientation\footnote{Thus our objects correspond to isomorphism classes of compact \emph{oriented} 1-dimensional manifolds with $\partial$.}. A morphism between two such objects is represented by a Riemann surface with boundary, with an embedding of the two objects into its boundary. We now make these ideas more precise.

\vspace{2mm}

We denote by $\hh, \overline \hh$ the open and closed upper half-planes, respectively. A map $f : U \to \cc$ defined on a relatively open subset $U \subseteq \overline \hh$ will be called \textbf{holomorphic} if it is continuous and holomorphic on $U \cap \hh$.

For us, \textbf{a Riemann surface} $X$ is a smooth oriented manifold with boundary with an atlas of local diffeomorphisms onto open subsets of $\overline \hh$ (or sometimes $\cc$, to ease notation) such that the transition maps are bi-holomorphic.

\begin{remark} The following somewhat technical observation will be useful:
one could start with a Hausdorff topological space\footnote{2nd countability follows from the other assumptions in this case.} $X$, and consider atlases of \emph{homeomorphisms} $f_\alpha : U_\alpha \to V_\alpha$ onto open subsets of $\overline \hh$ such that the transition maps are biholomorphic in the sense above.

Under these conditions we want to claim that the $f_\alpha$'s can be used to define a unique smooth structure of a manifold with boundary on $X$ such that all holomorphic maps are smooth.

Clearly it is enough to show that any transition map $f_2 \circ f_1^{-1}$ is smooth at $f_1(p)$, where $p \in U_1 \cap U_2$ is some boundary point. Restricting as necessary we may assume $U_1 = U_2 =: U$ is a simply connected domain.

The transition map $f_2 \circ f^{-1}_{1|f_1(U) \cap \hh}$ is then a Riemann mapping between the domains $f_i(U) \cap \hh$. Note that the real boundary segments $f_i(\partial X \cap U)$ are free boundary segments\footnote{Cf. \textsection \ref{subsec:reflection} for the definition of a free boundary segment. See also \citep[pg. 233]{ahlforsCA} where this is called a ``free one-sided boundary segment''.}. The reflection principle now implies that $f_2 \circ f^{-1}_{1|f_1(U) \cap \hh}$ extends holomorphically to some neighbourhood of these real segments. Since the continuation to the boundary is unique the same is true for $f_2 \circ f_1^{-1}$. In particular this function is smooth at $f_1(p)$.

In fact, we see that all sensible definitions for a Riemann surface with boundary are equivalent. We can require the transitions to be: continuous up to boundary / smooth up to boundary / holomorphic across the boundary...
\end{remark}

\begin{definition}
\label{def:OCR}
\textbf{An open-closed Riemann surface $X$ from $(C_-,O_-)$ to $(C_+,O_+)$}, denoted also $X : (C_-,O_-) \Rightarrow (C_+,O_+)$, consists of the following data:
\begin{enumerate}
\item A Riemann surface with boundary.
\item Smooth embeddings $\alpha_- : (S^1)^{C_-} \hookrightarrow \partial X$ and $\alpha_+ : (S^1)^{C_+} \hookrightarrow \partial X$.
\item Smooth embeddings $\beta_- : [0,1]^{O_-} \hookrightarrow \partial X$ and $\beta_+ : [0,1]^{O_+} \hookrightarrow \partial X$
\end{enumerate}
we require that $X$ be compact; that the images of the parametrization maps are disjoint; and that $\alpha_+, \beta_+$ preserve orientation (with respect to the induced outward normal orientation on the boundary) and $\alpha_-, \beta_-$ reverse orientation.
\end{definition}

Two OCR-surfaces $X,Y$ are considered equivalent if there's a biholomorphism $f : X \to Y$ which respects the parameterizations of the boundary. We can define $M_\Lambda((C_+,O_+),(C_-,O_-))$ as the set of equivalences classes\footnote{Some care needs to be taken if one wants to show that this is a proper set. At any rate this definition will not be used; in what follows we will give a much more tangible construction of this set.}.

\vspace{3mm}
We will sometimes want to forget some of the structure of an OC-Riemann surface, e.g. talk about \textbf{OC-smooth surfaces} and the like. At any rate, we will only consider oriented surfaces and orientation-preserving morphisms. For a definition of orientation in the category of topological surfaces see \citep[pg. 8-9]{LV}.

An (orientation-preserving) mapping between surfaces with parametrized boundary is said to be \textbf{rel-$\partial$} if it respects the parametrization of the boundary.

\chapter{Smooth Deforming and Stretching}
\label{chapt3}
\label{sec:smoothdef}

To proceed with the definition of Segal's category we need to define a topology on the set of morphisms, as well as continuous composition maps. Clearly the representation of a surface by a collection of charts is unwieldy and redundant. What is needed is a way to bundle OCR-surfaces together.

Let us describe how one might do this by thinking of OCR-surfaces as smooth deformations of one another. This construction will allow us to topologize the morphism space, and will be analogous to the construction we'd utilize eventually. It will also exhibit clearly why smooth maps are deficient in defining the sewing of Riemann surfaces, and thus motivate the introduction of quasiconformal mappings in the next chapter.

We say that two OCR-surfaces $(C_-,O_-) \Rightarrow (C_+, O_+)$ have \textbf{the same OC-smooth type} if there's a diffeomorphism rel-$\partial$ between them. We sometimes record the in and out objects when discussing OC-smooth types, e.g. $\tau : (C_-,O_-) \Rightarrow (C_+, O_+)$. For the classification of OC-smooth types in terms of discrete data, see corollary \ref{cor:type-data}.

Let us focus on one OC-smooth type $\tau$, and fix an OC-smooth model $\Sigma = \Sigma_\tau$ for it. Let $\mathcal{J}(\Sigma_\tau) \to \Sigma_\tau$ be the subbundle of $\mbox{End}(T\Sigma) \to \Sigma_\tau$ whose fiber over the point $p \in \Sigma_\tau$ is

\[\mathcal{J}(\Sigma_\tau)_p = \{J \in \mbox{End}(T\Sigma)_p | J^2 = -\mbox{Id}\}\]

We call smooth sections $C^\infty(\Sigma,\jj(\Sigma_\tau))$ \textbf{almost complex structures on $\Sigma_\tau$}. It is easy to see that every conformal structure on $\Sigma$ which agrees with the smooth structure gives rise to an almost-complex structure corresponding to ``multiplication by $i$'' in the tangent space, and that this almost complex structure determines the complex structure completely (i.e. this map is injective). It is a non-trivial fact that this map is also surjective - every smooth almost-complex structure on $\Sigma$ is obtained from a complex structure in this way. See \citep{chern} for a proof. Theorem \ref{thm:integration} is the quasiconformal analogue of this result which we will use.

Let $X$ be an OCR-surface of type $\tau$, and $f : \Sigma \to X$ a diffeomorphism rel-$\partial$. Complex multiplication by $i$ gives a smooth section $J_X\in C^\infty(X, (\mathcal{J}(X)))$, and by pulling back we obtain an almost complex structure on $\Sigma$, $f^* J_X \in C^\infty(\Sigma, (\mathcal{J}(\Sigma)))$.

Clearly, if $f_i : \Sigma \to X_i$, $i = 1,2$ satisfy $f_1^* J_{X_1} = f_2^* J_{X_2}$, then $f_2 \circ f_1^{-1} : X_1 \to X_2$ is a biholomorphism, though not all biholomorphisms between $X_1$ and $X_2$ are obtained in this way. $X_1$ and $X_2$ are biholomorphic if and only if there exists some diffeomorphism rel-$\partial$ $\psi : \Sigma \to \Sigma$ such that $f_1^* J_{X_1} = \psi^* f_2^* J_{X_2}$. Thus we are lead to consider the quotient
\[\mm_{\Lambda}(\tau) := C^\infty(\Sigma, (\mathcal{J}(\Sigma))) / \mbox{Diffeo}(\Sigma),\] where $\mbox{Diffeo}(\Sigma)$ denotes the group of self-diffeomorphisms of $\Sigma$ rel-$\partial$. The above discussion shows that the points of type $\tau$ of the morphism space are in bijective correspondence with $\mm_\Lambda(\tau)$. The total morphism space is then the disjoint union over OC-smooth types:

\[\mm_\Lambda((C_-,O_-),(C_+,O_+)) := \coprod_{\{\tau | \tau : (C_-,O_-) \Rightarrow (C_+,O_+)\}} \mm_\Lambda(\tau)\]

Note that this set has a natural topology, induced from the $L^\infty$ operator norm  on sections of $\mbox{End}(T\Sigma)$ through the canonical constructions of subset, quotient and disjoint union topologies.

Let $\tau_1 : (C_-,O_-) \Rightarrow (C_0,O_0)$ and $\tau_2 : (C_0,O_0) \Rightarrow (C_+,O_+)$ be composable OC-smooth types, and let us try to define the continuous composition maps $\mm_\Lambda(\tau_1) \times \mm_\Lambda(\tau_2) \to \mm_\Lambda(\tau_2 \circ \tau_1)$ (cf. remark \ref{rem:composing-types} about the composition of OC-smooth types). We focus first on the case that $O_0 = 0$, i.e. the sewing is along closed boundary circles only.

If we fix smooth models $\Sigma_1,\Sigma_2,\Sigma_{12}$ (e.g., by fixing $\Sigma_{12}$ and then drawing closed curves that bisect it) and then try to glue two almost complex structures together, we immediately run into trouble since the sections may not even have the same limit near the cut. Note, however, that there are many smooth structures on the topological union $\Sigma_1 \cup \Sigma_2$ that are compatible\footnote{We say a smooth structure on $\Sigma_1 \cup \Sigma_2$ is ``compatible'' with the smooth structures on $\Sigma_1$ and $\Sigma_2$ if the restrictions of every smooth map on $\Sigma_1 \cup \Sigma_2$ to $\Sigma_1$ and to $\Sigma_2$ are smooth.} with the smooth structures on the pieces (to see this, consider first the various smooth structures on $[-1,1]$ which are compatible with the usual smooth structures on $[-1,0]$ and $[0,1]$). In fact, it is not hard to see that there are smooth structures on $\Sigma_1 \cup \Sigma_2$ for which the almost complex structures can be glued at least \emph{continuously} - these are precisely the smooth structures which identify the vector field $\frac{d\alpha_1}{dt}$ with $\frac{d\alpha_2}{dt}$ and the vector field $\left(f_1^*J_{X_1}\right) \frac{d\alpha_1}{dt}$ with $\left(f_2^*J_{X_2}\right) \frac{d\alpha_2}{dt}$. Here $\frac{d\alpha_i}{dt}$ is the unit velocity vector field along $\partial X_i$ defined by the parametrization $\alpha_i$.

So it seems plausible, and in fact in appendix \ref{sec:smooth sewing} we show it to be true, that there is a smooth structure on $\Sigma_1 \cup \Sigma_2$ compatible with the smooth structures on $\Sigma_1,\Sigma_2$, which allows the two almost complex structures on the pieces to be glued \emph{smoothly}, and thus define an almost complex structure on the sewed surface. But even considering only the continuity requirement it is evident that for different pairs of almost complex structures, different smooth structures will be needed to facilitate the gluing. This implies that we need to step out of the smooth category in order to be able to bundle all of the relevant smooth structures together, and quasiconformal mappings will allow us to do that. Another intimately related motivation for turning to the quasiconformal category is that the gluing of smooth maps is not necessarily smooth, so if we want the gluing of almost complex structures to descend to the quotient, i.e. to the morphism space, we must look for a category in which the gluing of two automorphisms is again an automorphism.

\begin{figure}[hb]
  \centering
\includegraphics[width = 30mm]{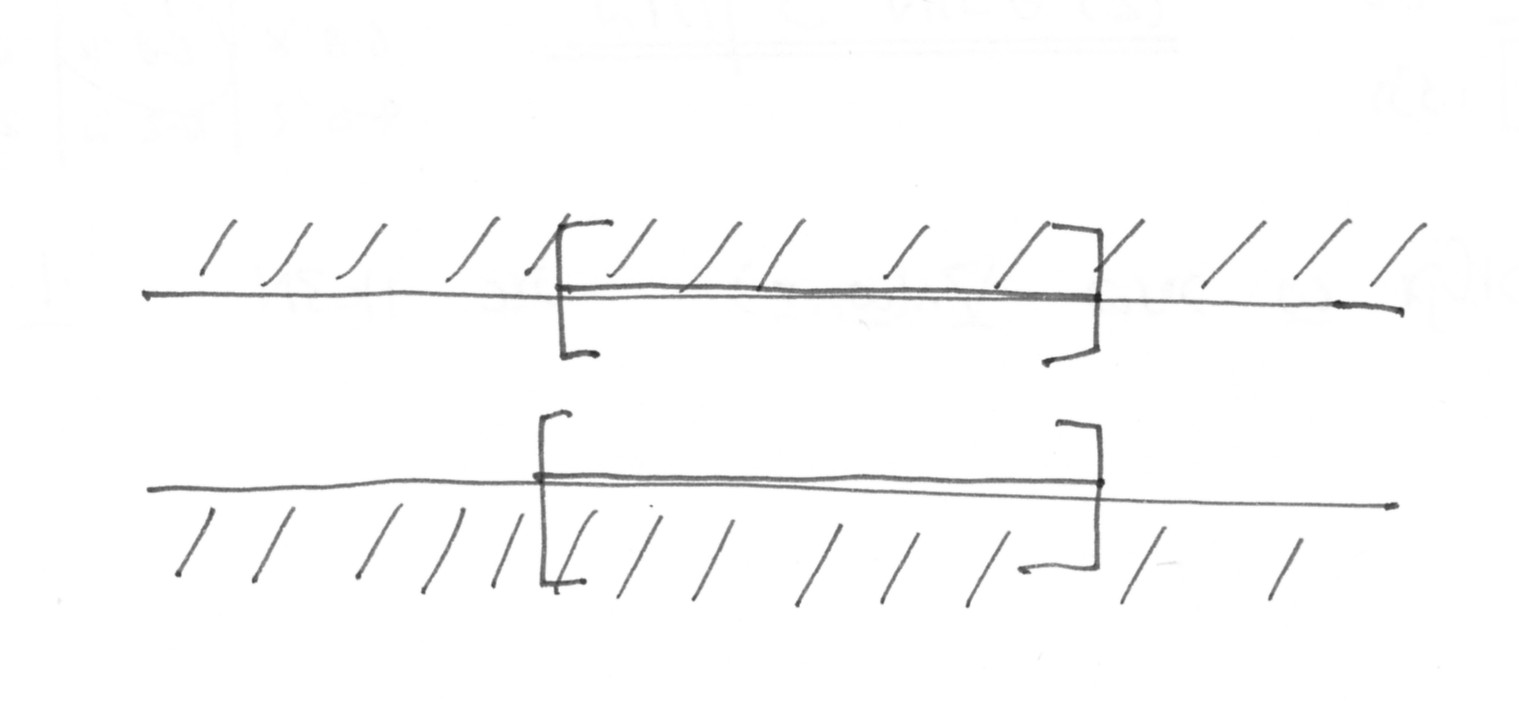}
\includegraphics[width = 30mm]{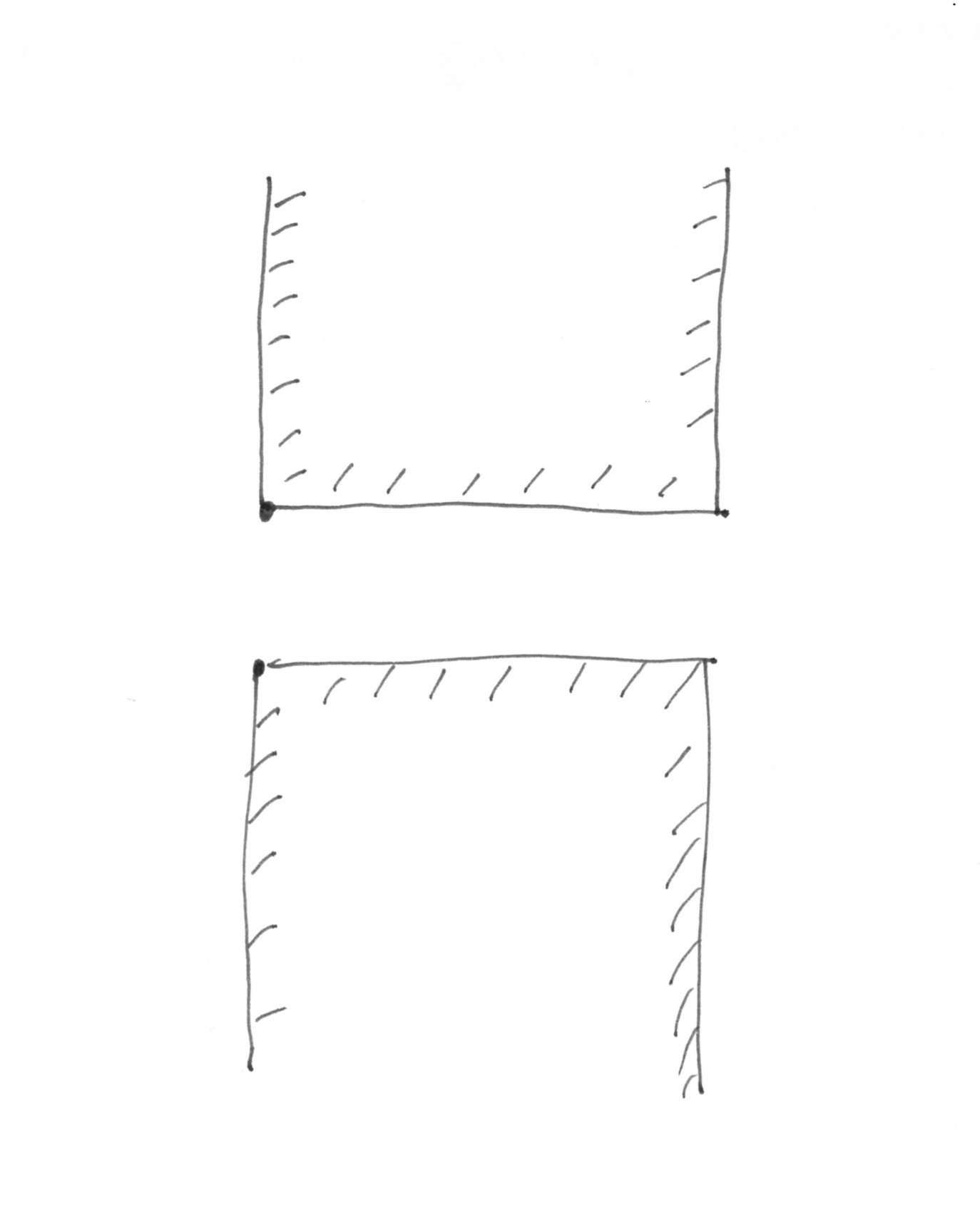}
\includegraphics[width = 30mm]{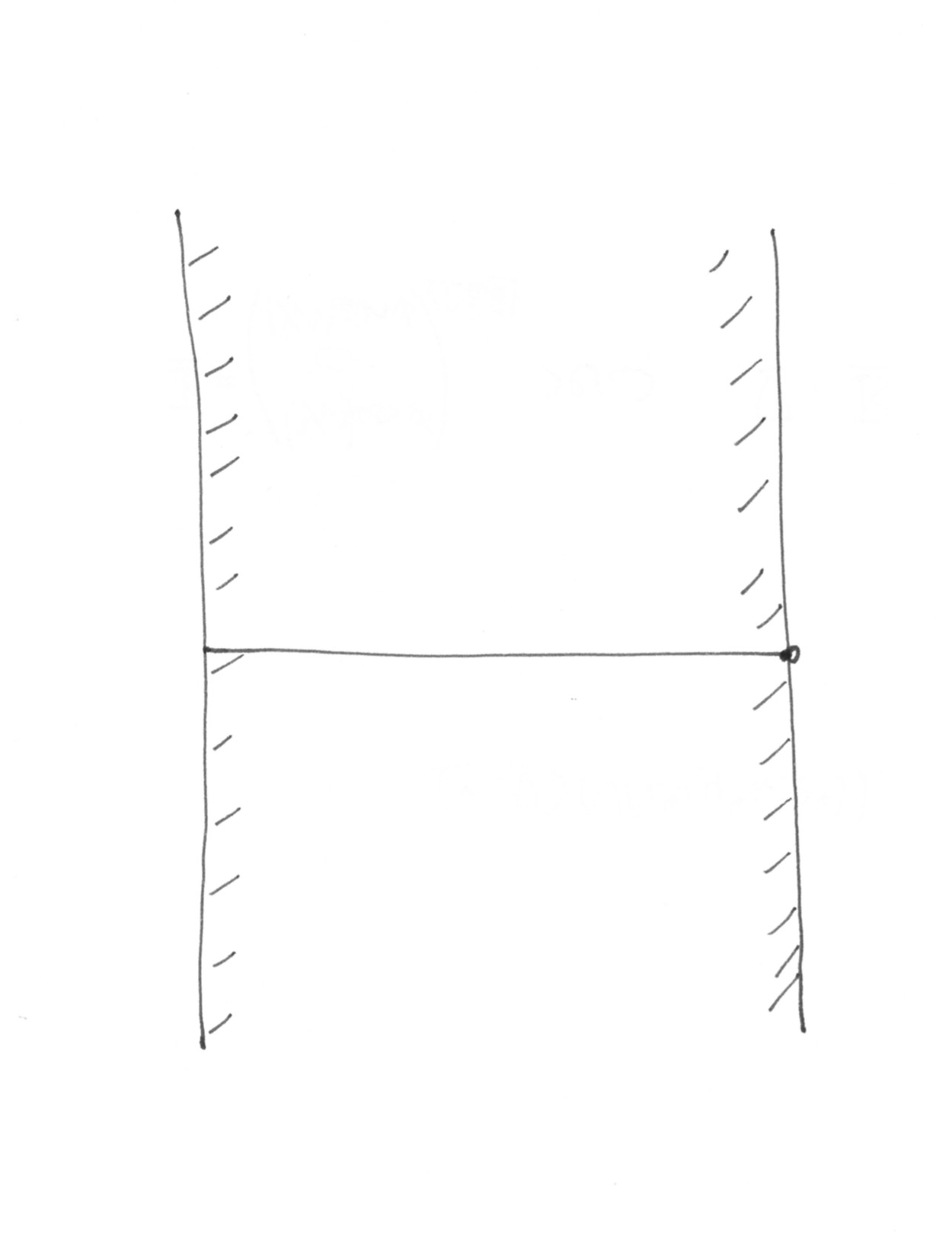}
  \caption[corners]
  {When gluing along open boundary segments we must first introduce corners.}
\label{fig:corners}
\end{figure}

Our final motivation for turning to QC mappings is related to the open boundary components. While it is intuitively - and topologically - quite clear what we mean by sewing two surfaces along an open boundary component, there is in fact \emph{no} smooth structure on $\Sigma_{12}$ which is compatible with $\Sigma_1$ and $\Sigma_2$ in this case. This is because we need to introduce corners in the smooth boundaries of $\Sigma_1$ and $\Sigma_2$ before we can piece them together to form a new smooth boundary (see figure \ref{fig:corners}). More formally, if we consider a chart of the glued surface and the two pieces near an endpoint of an open boundary segment we see that the transitions cannot be diffeomorphic - if we look at each surface in the separated picture, for every point $p$ (including boundary points) there's a smooth path whose midpoint is $p$. This is not true if we look at any one piece of the glued picture - there's no smooth path whose midpoint is the endpoint of the boundary segment. In contrast to this, quasiconformal maps can introduce corners, and so the QC category facilitates sewing along open boundaries as well.

\begin{figure}[hb]
  \centering
  \includegraphics[width=40mm]{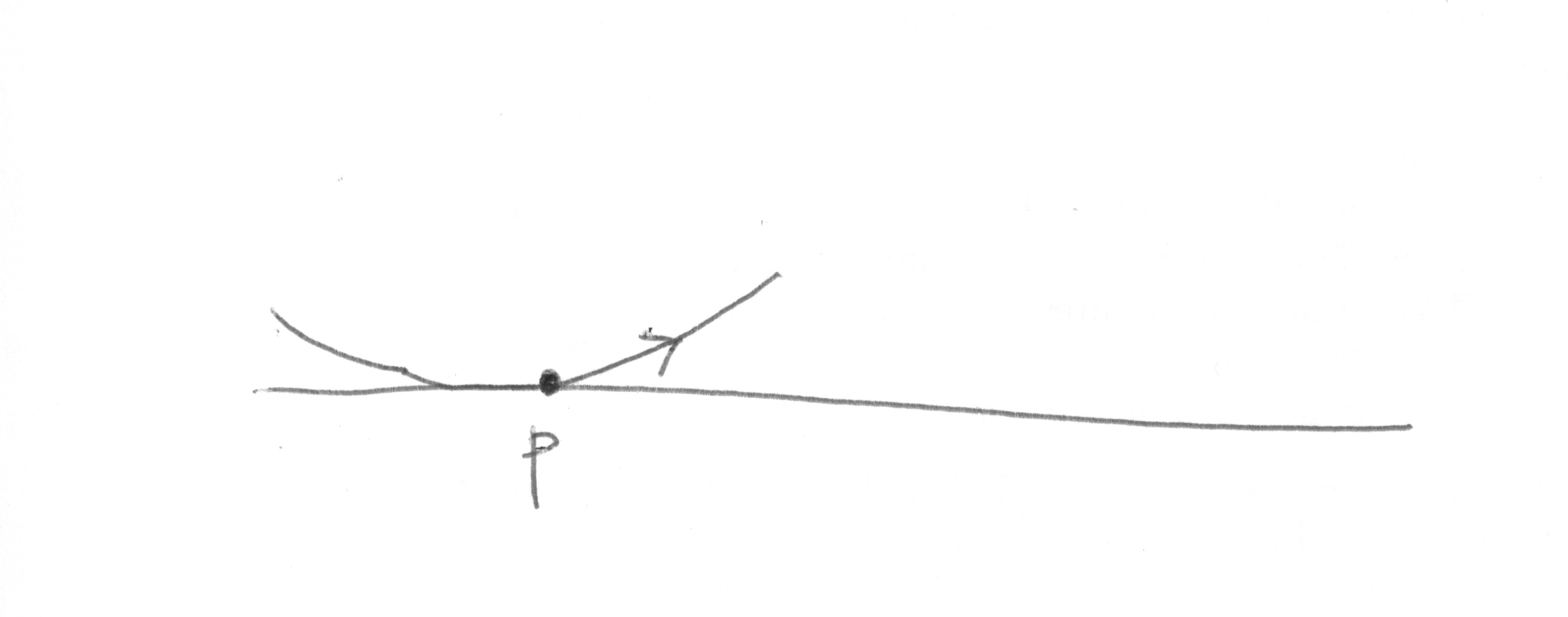}
  \includegraphics[width=40mm]{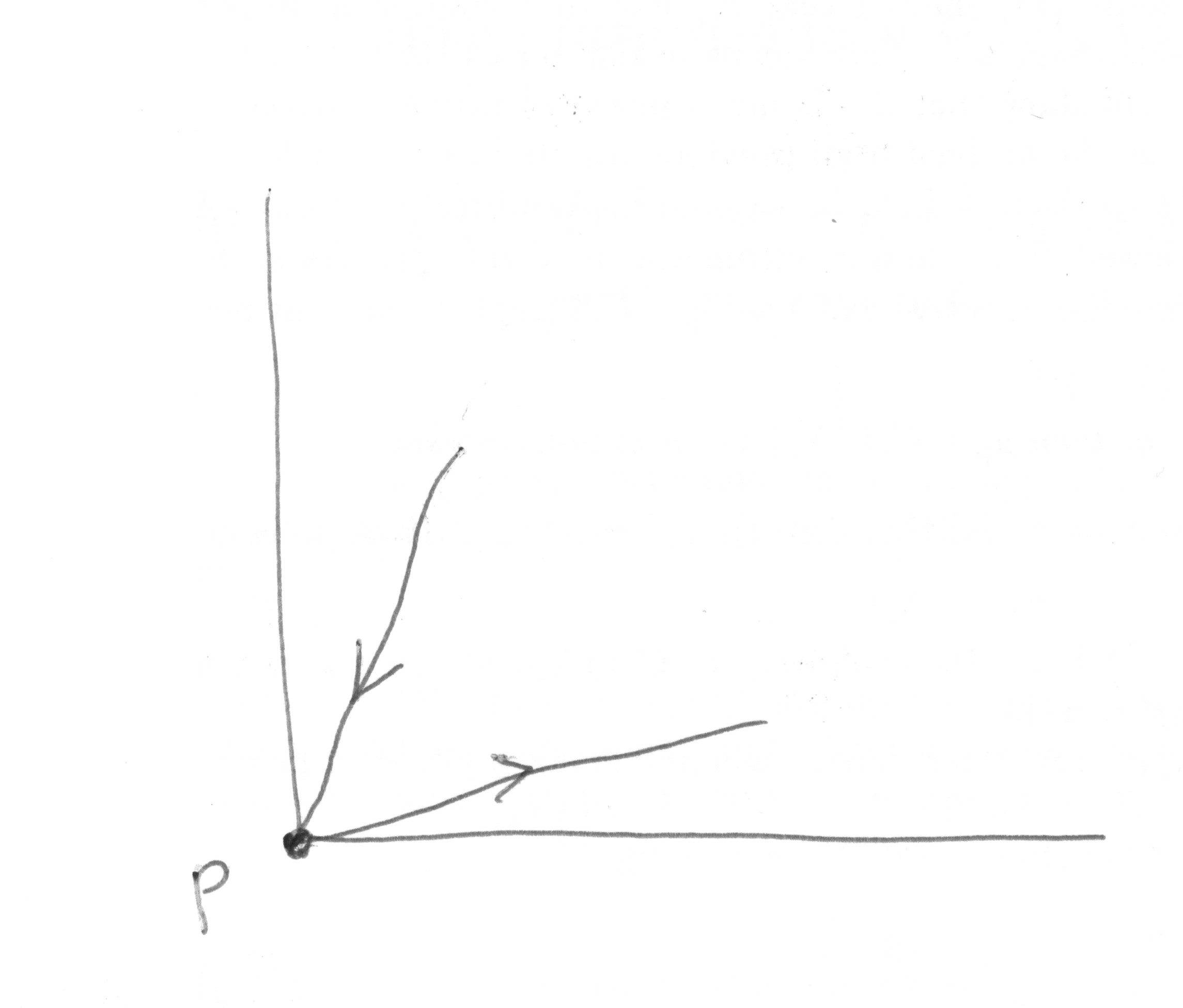}
  \caption[Path Midpoints]%
  {Before sewing, every point $p$ is a midpoint of a path. The same cannot be said for the corners after the sewing.}
\end{figure}

\chapter{Quasiconformal Mappings}
\label{chapt4}
\label{sec:QCmappings}

In this chapter we highlight the necessary ingredients of the theory of quasi-conformal mappings, and prove a few lemmas which demonstrate how they may be utilized.

\citep{LV} is a very thorough and well-written introduction to the theory of quasiconformal mappings. In what follows, we will continuously refer to this book simply as [LV].

\defcitealias{LV}{LV}

The lecture notes \citep{ahlfors-lectures} can be used as a more rapid introduction both to QC mappings and the basics of Teichm\"uller theory.

\section{The Geometric Definition}
Let $G,G'$ be open domains in the plane. Here and in what follows, when we speak of ``the plane'' we always mean the extended plane, i.e. the Riemann sphere. A mapping for us will always be an orientation-preserving homeomorphism $G \to G'$, unless otherwise stated. Roughly speaking, a mapping is quasiconformal if it is close to being conformal, and hence does not distort angles too much. To make this more precise we need the notion of a topological quadrilateral.

\begin{definition}
A \textbf{quadrilateral} is a Jordan domain\footnote{A Jordan curve is a subset of the plane which is homeomorphic to the circle. A Jordan domain is a domain (connected, open subset of the plane) whose boundary is a Jordan curve. By the celebrated Jordan theorem, the complement of a Jordan curve is disconnected and consists of two Jordan domains.} $Q$ with a 4-tuple of distinct points $(z_0,z_1,z_2,z_3)$ on its boundary. We assume that these points (or ``vertices'') are ordered with positive orientation relative to $Q$\footnote{More precisely, choose a conformal mapping $Q \to D$ where $D$ is the unit disc. The images of the points should appear in cyclic order when we traverse the circle in a counter-clockwise fashion.}. We will speak of the quadrilateral $Q(z_0,z_1,z_2,z_3)$ or sometimes just $Q$, for short.

A \textbf{rectangle} is a quadrilateral whose sides are parallel to the coordinate axes and whose vertices coincide with the corners.
\end{definition}

By the Riemann mapping theorem every quadrilateral can be mapped conformally onto any other Jordan region in the plane, and this map extends continuously to a homeomorphism of the boundaries\footnote{\label{footnote:OTC} For the second assertion, see \citep{osgood-taylor-caratheodory}. In fact the proof generalizes to QC mappings, see \citepalias[pg. 81]{LV}}. If $\phi$ is such an extended mapping between $Q$ and $\phi(Q)$, We consider the quadrilateral $Q(z_0, z_1,z_2,z_3)$ equivalent to $\phi(Q)(\phi(z_0), \phi(z_1), \phi(z_2), \phi(z_3))$.

\begin{definition}
\textbf{The module} of $Q$ is the unique positive number $M$ such that the rectangle $R(0, M, M + i, i)$ is equivalent to $Q$. We will also write $M(Q)$ for the module of the quadrilateral $Q$.
\end{definition}

To see that this is well-defined, let $Q$ be any quadrilateral with distinct marked points on the boundary. We can map $Q$ to the upper half plane, and then use a real M\"obius transformation to send $Q(z_1,z_2,z_3,z_4)$ to $Q(-1,0,1,x)$. Finally, we apply the Schwarz-Christoffel formula\footnote{See \citep[pg. 236]{ahlforsCA}.} to map the upper half $z$-plane to a rectangle of the form $R(0, M, M + i, i)$ in the $w$-plane:
\[w = C\int_0^z \frac{dz'}{\sqrt{z'(z'^2-1)(z'-x)}} + C'.\]
Clearly, the equivalence class of $Q$ then also contains all similar rectangles.

We still need to show uniqueness. Using the reflection principle, extend any mapping between rectangles $R$ and $R'$ to a conformal mapping of $\cc$; such a mapping must be of the form $z \mapsto az + b$ and we see that equivalence classes cannot contain more than one similarity class of rectangles. In particular, every equivalence class contains a unique rectangle of the form $R(0,M,M+i,i)$. \qed

\vspace{3mm}
We say that a quadrilateral is in $G$ if its closure is contained in $G$.

\begin{definition}
Let $K$ be a finite, positive number. A mapping $f : G \to G'$ is called $K$-quasiconformal if for any quadrilateral $Q$ in $G$ we have
\[\frac{1}{K}M(Q) \leq M(Q') \leq K M(Q)\]
where $Q'$ is the quadrilateral in $G'$ which is the image of $Q$ under $f$.
\end{definition}

We will sometimes find it sufficient to say that $f$ is quasiconformal when it is $K$-quasiconformal for some finite $K$. The minimal $K$ for which $f$ is $K$-QC is called the maximal dilatation of $f$.

We list a few immediate consequences of the definition
\begin{lemma}
\label{lemma:QCproperties}

\begin{enumerate}
\item $f$ and $f^{-1}$ have the same maximal dilatation.
\item The class of $K$ q.c. mappings is invariant under conformal mappings.
\item The composite of $K_1$ and a $K_2$ QC mapping is $K_1 K_2$ QC.
\end{enumerate}
\end{lemma}

It is also possible to show that a map which is $1$-QC is in fact conformal (the converse is immediate).

\section{The Analytic Definition}
\label{subsec:analyticdef}
Quasiconformal mappings turn out to be very well-behaved analytically. In this chapter we will give necessary and sufficient analytic conditions for a function to be quasiconformal.

Let $D,D'$ be two domains in $\cc$, and let $f : D \to D'$ be an orientation-preserving mapping which is regular (differentiable with nonzero jacobian) at $p \in D$. Suppose we try to pull back the usual complex structure on $D'$ to $D$ along $f$. We'd like to find an economic way to describe that part of $df_p$ which determines the almost complex structure at $p$. We can always write $df_p(z) = A z + B \overline z$, where $A = A(p), B = B(p)$ are complex numbers. It is then an easy matter to check that the complex structure is determined uniquely by the \textbf{complex dilatation} $\mu_f(p) = \frac{B}{A}$. This is also written as $\mu = \frac{f_{\overline z}}{f_z} =\frac{df/d\overline z}{df / dz}$\footnote{The notation $A = \frac{df}{dz}, B = \frac{df}{d\overline z}$ does not refer to a limiting process.}. Thus when we compare almost complex structures, it is more natural to consider the complex dilatation rather than the tensor $J \in \End(TX)$ which we encountered in chapter \ref{sec:smoothdef}. For orientation preserving diffeomorphisms we have $|A| > |B|$, so $|\mu| < 1$.

One calculates that if $f$ is a real affine mapping, then $\mu$ is constant in $D$ and is related to the maximal dilatation $K$ of $f$ by $|\mu| = \frac{K - 1}{K + 1}$. Thus we can think of the absolute value of $\mu$ as an infinitesimal measure of the dilatation of $f$. For more precise results along these lines, see \citepalias[pg. 197]{LV}.

\begin{definition}
A mapping $f : G \to G'$ is called \textbf{absolutely continuous on lines} if for every rectangle $R$ in $G - \{\infty\} - \{f^{-1}(\infty)\}$, $f$ is absolutely continuous on almost every horizontal and almost every vertical line segment of $R$.
\end{definition}

\begin{lemma}
Functions which are ACL are differentiable almost everywhere.
\end{lemma}
\begin{proof}
Note that an absolutely continuous function $[a,b] \to \cc$ is differentiable almost everywhere\footnote{This is proved in most textbooks on measure theory, and also in \citepalias[III 2.7]{LV}.}. Then the following striking result of Gehring and Lehto\footnote{For the proof, see \citepalias[III, Theorem 3.1]{LV}.} establishes the claim.

\emph{Let $G$ and $G'$ be domains in the finite plane and $f : G \to G'$ a homeomorphism having finite partial derivatives almost everywhere in $G$. Then $f$ is differentiable almost everywhere in $G$.}\end{proof}

\begin{theorem}
\label{thm:analytic-geometric}
An orientation-preserving homeomorphism $f : G \to G'$ is $K$-quasiconformal if and only if:

\begin{enumerate}[(a)]
\item $f$ is absolutely continuous on lines and
\item $|\frac{df}{d\overline z}| \leq \frac{K-1}{K+1}|\frac{df}{dz}|$ almost everywhere ($\frac{df}{d\overline z}$ and $\frac{df}{dz}$ are defined where $f$ is differentiable, which is a.e. given (a)).
\end{enumerate}
\end{theorem}

For the proof, we refer the reader to \citepalias[pg. 168]{LV}.

Note that while the analytic definition is easier to work with, some of the so-called immediate consequences of the geometric definition are not so easy to prove. Usually the trouble is in proving absolute continuity on lines.

On the other hand, there are results which are easier to prove using the analytic definition. The following useful lemma says that $K$-quasiconformality is a local property. This is proved in \citepalias[pg. 48]{LV} using geometric estimates, but the analytic definition allows us to supply an easier proof.

\begin{lemma}
\label{lemma:localglobal}
An orientation-preserving homeomorphism $w$ of a domain $G$ is $K$-quasiconformal iff it is locally $K$-quasiconformal. That is, if for every point $z \in G$ there exists some open $U$, $z \in U \subset G$ s.t. $w_{|U}$ is $K$-quasiconformal.
\end{lemma}

\begin{proof}
If $f$ is QC then it is locally QC.

For the other direction note that a map which is locally $K$-QC is differentiable a.e. and its complex dilatation is bounded by $\frac{K-1}{K+1}$. It remains to show that it is ACL. For this let $R = [a,b] \times [c,d]$ be a rectangle in $G - \{\infty\} - \{w^{-1}(\infty)\}$. By standard Lebesgue-number arguments we can choose partitions $a = x_0 < x_1 < \cdots < x_n = b$ and $c = y_0 < y_1 < \cdots < y_n = d$ so fine that $w$ is $K$-QC in an open subset containing each $R_{ij} := [x_{i-1},x_i] \times [y_{j-1},y_j]$, $i,j = 1,...,n$.

A map $\rr \to \cc$ which is absolutely continuous on $[x_0,x_1]$ and on $[x_1,x_2]$ is absolutely continuous on $[x_0,x_2]$, so the set
\[\{y \in [y_{j-1},y_j] | w_{|[a,b] \times \{y\}} \mbox{ is absolutely continuous}\}\]
contains all points of $[y_{j-1}, y_j]$ perhaps except for points in the union of $n$ sets of measure 0. We conclude that the union of these sets for $j = 1,...,n$ differs from $[c,d]$ in a set of measure 0. The vertical segments can be dealt with symmetrically and the result is proven.
\end{proof}

\section{The Reflection Principle}
\label{subsec:reflection}

We will want to consider QC mappings for subsets of the plane which are not necessarily open.

A homeomorphism $w : S \to S'$ between two arbitrary sets in the plane will be called \textbf{$K$-quasiconformal} if it admits a $K$-quasiconformal extension to an open neighbourhood of every point in the domain. Note that if $S$ happens to be open then $S'$ is open too\footnote{See theorem V.13.2 in \citep{newman}.} and in this case this definition is equivalent to the previous one by lemma \ref{lemma:localglobal}.

While this definition has the merit of being a local property, it is somewhat awkward to work with. Fortunately we will not need to use it in full generality, as we're mostly interested in the case when $S,S'$ are relatively open subsets of $\overline \hh$. In this case we will to show that a map is $K$-quasiconformal (locally) iff it admits a $K$-QC extension to open neighbourhoods of $S,S'$ in the plane. This calls for an application of the reflection principle, which can be generalized to QC mappings as follows.

Recall that reflection in a circle $C$ (or a circular arc representing it) is an involution of the plane defined by $z \mapsto C(z) := T^{-1}(\overline {T(z)})$ where $T$ is any M\"obius transformation which sends $C$ to the real axis.

Let $G$ be a domain in the plane. An open Jordan arc (or curve) $C$ on the boundary of a domain $G$ is called a \textbf{free boundary arc of $G$} if:
\begin{enumerate}
\item Every point on $C$ contains an open disc whose intersection with the boundary of $G$ is contained in $C$.
\item If $E$ is the connected component of the boundary of $G$ which contains $C$, then $E - C$ is connected.
\end{enumerate}

We say that $G$ admits reflection in an open circular arc $C$ if (a) $C$ is a free boundary arc of $G$, and (b) $G$ does not intersect $C(G)$, the domain which is obtained from $G$ by reflection in $C$. In this case $G \cup C \cup C(G)$ is a domain.

It is important to note that if $U$ is any relatively open subset of $\overline \hh$ then $U \cap \hh$ admits reflection in any connected component of $U \cap \rr$.

\begin{theorem}[The Reflection Principle]
\label{thm:reflection}
Let $G$ and $G'$ be two Jordan domains which admit reflection in $C$ and $C'$, respectively. If $w : G \to G'$ is $K$-QC mapping which extends to a homeomorphism of $G \cup C$ onto $G \cup C'$, then $w$ extends to a $K$-QC mapping $\hat w : G \cup C \cup G^* \to G' \cup C' \cup G'^*$ according to the rule $w(C(z)) = C'(w(z))$.
\end{theorem}

\begin{remark}
It immediately follows from our assumptions that $\hat w$ is a homeomorphism which is differentiable a.e. and whose complex dilatations satisfy the required bound a.e. in $G$ and in $G'$, and hence a.e. in $G \cup C \cup G'$.

This does not establish that $\hat w$ is $K$-QC, as we have not proved that the map is ACL. This gap seems to be non-trivial, even when $C$ is a real line segment, since it is not necessarily true that a function which is ACL in the upper and lower open half planes is ACL in the plane. To see this, consider the map $(x,y) \mapsto y \sin{1/y}$, which is absolutely continuous on every vertical interval $\{x\} \times [-1, -a]$ and $\{x\} \times [a,1]$ for $a > 0$, but is not absolutely continuous on $[0,1] \times [-1,1]$.

A proof of theorem \ref{thm:reflection} which relies on geometric estimates can be found in \citepalias[pg. 46]{LV}. The result there is stated under the assumption that $C$ and $C'$ merely ``correspond to one another under $w$''. This is easily shown to follow from our assumptions.
\end{remark}

Consider again the special case where $U,U'$ are open subsets of $\overline \hh$, and let $w : U \to U'$ be a homeomorphism. $w$ maps $U \cap \rr$ onto $U' \cap \rr$, since the real points are precisely those points whose removal from $U$ (resp., $U'$) does not alter the fundamental group. We may therefore apply the reflection principle to $w_{|U \cap \hh}$ to obtain

\begin{lemma}
\label{lemma:reflect-localglobal}
For a homeomorphism $w : U \to U'$, $U,U'$ relatively open subsets of $\overline \hh$, the following conditions are equivalent:

\begin{enumerate}[(a)]
\item $w$ is $K$-QC (locally around each point)
\item $w_{|U \cap \hh}$ is $K$-QC
\item $w$ admits a $K$-QC extension to an open neighbourhood of $U$ in the plane.
\end{enumerate}
\end{lemma}

\section{Bending and Gluing}
\label{subsec:bending and gluing}
We now show that quasiconformal mappings overcome two obstacles that we stumbled upon in the previous chapter when we restricted ourselves to working with smooth maps. Namely, we want to show that QC mappings allow us to introduce corners in the boundary, and that they can be glued along straight lines.

\begin{lemma}
\label{lemma:corners}
The mapping $\sigma : \overline \hh \rightarrow Q$ given by $\sigma(z) = \sqrt{z}$ is quasiconformal. Here $Q$ is the closed first quadrant of the plane.
\end{lemma}

\begin{proof}
Let $\phi : S^1 \to S^1$ be an orientation-preserving $C^1$-homeomorphism with $\phi(e^{i\theta}) = e^{i\theta/2}$ for $0 \leq \theta \leq \pi$. We can extend this to a homeomorphism $\sigma : \cc \to \cc$ by setting $\sigma(r e^{i \theta}) = r^{1/2}\phi(e^{i\theta})$. $\sigma$ is a homeomorphism which is $C^1$ away from the origin. It is ACL, since it is differentiable on all vertical and horizontal line segments that do not pass through the origin. Its differential multiplies radial vectors by $\frac{1}{2\sqrt{r}}$ and multiplies vectors tangent to circles centered at the origin by a real scalar $\in [\min_{S^1} \frac{d\phi}{d\theta} \frac{\sqrt{r}}{r}, \max_{S^1} \frac{d\phi}{d\theta} \frac{\sqrt{r}}{r}]$. We deduce that its maximal dilatation is $\max\left \{\frac{1}{2} \max_{S^1} \frac{d\phi}{d\theta},2 / \min_{S^1} \frac{d\phi}{d\theta}\right\}$. It follows that $\sigma$ is QC on $\cc$ by theorem \ref{thm:analytic-geometric}, and $\sigma(\overline \hh) = Q$, as desired.
\end{proof}

Let $Q_1,Q_2$ be the first and second closed quadrants in the plane (so $\overline \hh = Q_1 \cup Q_2$). Let $L = Q_1 \cap Q_2$ be the non-negative imaginary axis. Let $U,V$ be open subsets of $\overline \hh$, and set $U_i := U \cap Q_i$, $i = 1,2$.

\begin{lemma}
\label{lemma:local pushout}
Assume $f : U \to V$ is a homeomorphism. Then $f$ is QC iff $f_i := f_{|U_i}$ are QC for $i=1,2$.
\end{lemma}

\begin{figure}[hb]
  \centering
\includegraphics[width = 100mm]{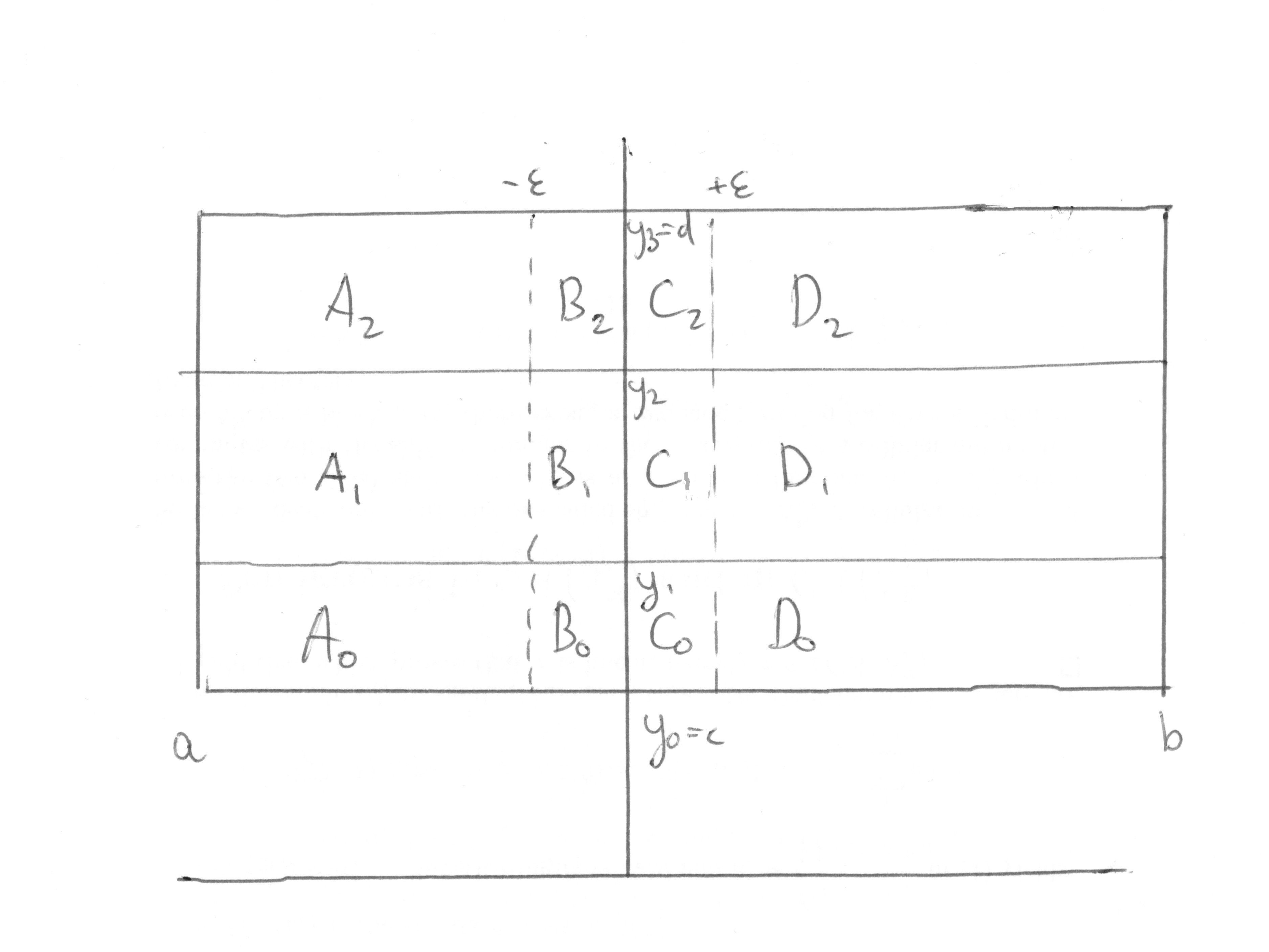}
\caption[refining]{By assumption, $f_1$ and $f_2$ are ACL on each of the depicted rectangles.}
  \label{fig:refining}
\end{figure}

\begin{proof}
If $f$ is QC then so are the $f_i$. For the other direction, assume that $f_{|U_i}$ are QC, and note that by lemma \ref{lemma:reflect-localglobal} it suffices to prove that $f$ is QC in $U \cap \hh$.

Clearly $f$ is differentiable almost everywhere in $U \cap \hh$ and $\|\mu_{f_1 \cup f_2}\|_\infty = \max(\|\mu_{f_1}\|_\infty,\|\mu_{f_2}\|_\infty)$. We must therefore show that $f$ is ACL.

Let $R := [a,b] \times [c,d] \subset U \cap \hh$ be a closed rectangle. Clearly we can focus our attention on the case when $a < 0 < b$. By compactness we can choose a partition $c = y_0 \leq y_1 \leq \cdots \leq y_n = d$ and $\epsilon > 0$ such that both $f_1$ and $f_2$ admit a $K$-QC extension to open neighbourhoods which contain $[-\epsilon,+\epsilon] \times [y_k, y_{k+1}]$ for every $k$. Let $A_k,B_k,C_k,D_k$ denote the rectangles obtained as the product of the vertical segment $[y_k,y_{k+1}]$ with the horizontal segments $[a,-\epsilon],[-\epsilon,0],[0,\epsilon]$ and $[\epsilon,b]$, respectively (see fig. \ref{fig:refining}). We find that $f_1 \cup f_2$ is absolutely continuous on a.e. horizontal and vertical segments of the $A_k,B_k,C_k,D_k$ and therefore, as in the proof of lemma \ref{lemma:localglobal}, $f_1 \cup f_2$ is absolutely continuous on a.e. horizontal and vertical segment of $R$.
\end{proof}

\section{Transformation Rule for Complex Dilatations}
\label{subsec:transformation}

Let $f$ and $g$ be quasiconformal mappings. A straightforward application of the chain rule gives:
\begin{equation}
\label{eqn:composition}
\mu_g \circ f = \frac{f_z}{\overline f_{\overline z}}\frac{\mu_{g \circ f} - \mu_f}{1 - \overline {\mu_f} \mu_{g \circ f} }\end{equation}
one should think of this as a transformation rule; we will use it to construct a geometric object in \textsection \ref{subsec:beltrami-diffs}. It is also good to note that this presents the relation between $\mu_g \circ f$ and $\mu_{g \circ f}$ as a M\"obius transfomration of the unit disc (which varies from point to point according to the derivatives of $f$). In particular, if $g_1,g_2$ are two QC functions then $d(\mu_{g_1},\mu_{g_2}) = d(\mu_{g_1 \circ f},\mu_{g_2 \circ f})$ where $d(s_1,s_2)$ is the essential supremum of the hyperbolic distance between $s_1(z)$ and $s_2(z)$.

\section{Prescribing Complex Dilatations}
We turn now to one of the deepest and most surprising result about QC mapping; indeed this may be thought of as ``the fundamental theorem of quasiconformal mappings''. There are several ways to prove this result, all of which require considerable effort. \citepalias[Chapter V]{LV} contains a complete proof as well as proof-sketches and references for other lines of reasoning.

\begin{theorem}
\label{thm:integration}
If $G$ is a domain and $\chi$ an arbitrary measurable function in $G$ with \[\sup_{z \in G} |\chi(z)| < 1,\]
then there exists a quasiconformal mapping $w$ of $G$ whose complex dilatation coincides with $\chi$ almost everywhere in $G$.
\end{theorem}

A few remarks are in order. First, this is often interpreted as an existence theorem for solutions of the Beltrami differential equation
\begin{equation}
\label{eqn:beltramiDE}
w_{\overline z} = \chi w_z.
\end{equation}
A solution in this context means either a ``solution almost everwhere'' or alternatively a solution ``in the weak sense'' (of distributions, i.e. functionals on smooth compactly supported test functions).

By the formula for composition of complex dilatations we see that if $w_1 : G \to G_1$ and $w_2 : G \to G_2$ are two maps with the same complex dilatation $\chi$ then $w_2 \circ w_1^{-1} : G_1 \to G_2$ has complex dilatation 0 a.e. in $G_1$. It follows that it is 1-QC, hence conformal. Thus the solution of equation (\ref{eqn:beltramiDE}) is unique up to suitable normalization. E.g., if $G$ is the entire (extended) plane then we may specify the images of $0,1$ and $\infty$.

\section{Sketch of Complex Sewing}
\label{subsec:sewing sketch}
The existence and uniqueness theorem has many applications. Let us give a rough sketch of how it may be used to sew complex structures. The precise details of this construction will be taken up in chapter \ref{sec:composition}.

Let $X_1 : (C_-,O_-) \Rightarrow (C_0, O_0)$ and $X_2 : (C_0, O_0) \Rightarrow (C_+,O_+)$ be two OCR-surfaces (we say that they are \textbf{composable}, since they represent composable morphisms). We can glue the two surfaces topologically according to the boundary parametrization. Let us denote the resultant surface by $X_{12}$ \footnote{We already mentioned the problems with giving this surface a smooth structure. In the next chapter we will see that we can in fact think of this as a QC surface. For our current purposes it is enough to consider it as a topological manifold.}. We claim that there's a unique conformal structure on $X_{12}$ such that for any open set $U \subset X_{12}$ a continuous function $U \to \cc$ is analytic iff $f_{|U_1}$ and $f_{|U_2}$ are analytic (according to the conformal structures on $X_1$ and $X_2$), where $U_i = X_i \cap U - \partial X_i$.

To define this structure we use smooth mappings (which, restricted to relatively compact sets, are quasiconformal) to align charts of $X_1$ and $X_2$ near the common boundary according to the parametrization, after introducing the appropriate corners for open boundary segments using lemma \ref{lemma:corners}. We keep track of the complex structures on $X_1$ and $X_2$ by the complex dilatation of the modified charts we are using. The complex dilatations may not agree on the seam (the common boundary), but since it has measure zero this does not pose a problem. We can now appeal to the existence theorem to find a conformal chart of the glued surface with the prescribed complex dilatations. This chart restricts to a conformal map on each piece away from the boundary - in fact it is determined uniquely (up to conformal transformation) by this requirement. This gives the desired result.

It is instructive to note that one can sew quite trivially when the reparametrizations $\alpha_2 \circ \alpha_1^{-1}$ and $\beta_2 \circ \beta_1^{-1}$ of circles and intervals of the common boundary are represented by analytic maps $\rr \to \rr$ in local holomorphic variables on $X_1$ and $X_2$; in this case we can simply juxtapose \emph{conformal} charts (since the aligning map is analytic) and thus there's no stretching to keep track of, hence no modified complex structure that needs to be integrated.

This is true even if we glue open boundaries since, as we remarked after the proof of lemma \ref{lemma:corners}, we can use the map $\sigma : z \mapsto \sqrt{z}$ to introduce the desired corners (in this case holomorphic maps on $X_{12}$ will not have analytic - or even smooth - continuations across the common boundary when restricted to the two pieces).

\chapter{Quasiconformal Surfaces}
\label{chapt5}
We want to define the objects and morphisms of the category of quasiconformal surfaces. We will then show how to bundle complex structures on quasiconformal surfaces in much the same way that we bundled complex structures on smooth surfaces in chapter \ref{sec:smoothdef}. Finally we will introduce the notion of a QC-structure with parametrization. This setup will facilitate the definition of the composition of morphisms in the next chapter.

\section{Quasiconformal Surfaces}
It follows from lemma \ref{lemma:QCproperties} that the set of QC mappings between open subsets of $\overline \hh$ is closed under compositions and taking inverses. This motivates the following definition.

\begin{definition}
Let $X$ be a 2nd countable compact Hausdorff topological space. A \textbf{quasi-conformal structure on $X$} is a maximal collection of homeomorphisms (``charts'') $X \supset U_i \overset{\varphi_i}\longrightarrow V_i$ such that

\begin{enumerate}
\item $V_i$ is an open set of $\overline \hh$.
\item the transition maps $\psi_{ij} = \varphi_j \circ \varphi^{-1}_{i| V_{ij}} :  V_{ij} \to V_{ji}$ are QC (cf. \textsection \ref{subsec:reflection}), where $V_{ij} := \varphi_i(U_i \cap U_j)$.
\item $X = \cup U_i$.
\end{enumerate}
Since the transition maps are orientation-preserving, the underlying surface is orientable. We include a choice of orientation as part of the data.
\end{definition}

Each transition map is $K$-quasiconformal for some $K$, but there is no universal bound on this $K$. Because $X$ is compact we may, however, choose a finite atlas (a finite set of charts whose domains cover $X$) and thus obtain a global bound on the $K$ of the transition maps. The assumption that $X$ is compact will be used time and again.

As usual, we call the preimages of the real axis \textbf{the boundary of $X$} and denote it by $\partial X$.

\begin{definition}
A map $f : X \to Y$ will be called a \textbf{QC isomorphism} (resp. \textbf{QC-embedding}) if it is a homeomorphism (resp., topological embedding) and
$\varphi_i \circ f \circ \psi_j^{-1}$ is QC (as a homeomorphism onto its image) for any pair of charts of $X$ and $Y$.

These are the only two kinds of \textbf{QC morphisms} we will consider. Note that if $f$ is an isomorphism then its inverse is indeed a QC morphism as well.
\end{definition}

It follows from the definition that a QC morphism is orientation-preserving. Also, although the definition calls for checking the QC condition for \emph{any} pair of charts, by transitivity it is in fact sufficient to check this for a pair of finite atlases of the domain and range.

\begin{lemma}
Isomorphisms are closed under compositions. Embeddings are closed under composition with isomorphisms and embeddings.
\end{lemma}

\begin{proof}
Let $f : X \to Y$ and $g : Y \to Z$ be isomorphisms or embeddings. Clearly $g \circ f$ is a topological isomorphism or embedding as the formulation of the lemma dictates.

It remains to show that if $\varphi_X : U_X \to V_X$ and $\varphi_Z : U_Z \to V_Z$ are charts of $X$ and $Z$ then there exists a constant $K$ such that $\varphi_Z \circ g \circ f \circ \varphi_X^{-1}$ extends to a  $K$-QC map about each point in $D := \varphi_X(f^{-1}(g^{-1}(U_Y)))$. Choose a finite atlas of $Y$, $\psi_i : U_i \to V_i$, $i = 1, \hdots,r$. We can set $K = \max_{i = 1,\hdots,r} K_g^i K_f^i$, where $K_g^i,K_f^i$ are bounds for the dilatations of $\varphi_Z \circ g \circ \psi_i^{-1}$ and $\psi_i \circ g \circ \varphi_X^{-1}$, respectively.
\end{proof}

We would like to relate QC structures to smooth structures. It turns out that every smooth structure determines a unique QC structure, as the next lemma shows.

A chart $\varphi : X \supset U \to V$ will be called \textbf{refined} if it is the restriction of a chart $\widetilde \varphi : \widetilde U \to \widetilde V$ with $\overline U \subset \widetilde U$.

\begin{lemma}
If $X$ is a smooth manifold with $\partial$ then the refined charts of $X$ belong to a unique QC structure on $X$. Moreover, if $f : X \to Y$ is a diffeomorphism or an immersion (of smooth manifolds with $\partial$) then it induces a QC-isomorphism or a QC-embedding, respectively.
\end{lemma}

\begin{proof}
If $f = \varphi_2 \circ \varphi_1^{-1}$ is the transition between smooth refined charts $\varphi_i$, then there is a positive constant $M$ such that $\frac{1}{M} \leq |\partial_\alpha f| \leq M$, where we set
\[\partial_\alpha f(z) := \lim_{r \to 0} \frac{f(z + r e^{i \alpha}) - f(z)}{r}\]
for the directional derivative of $f$. In terms of the directional derivatives the maximal dilatation is given by $K = \sup_z \frac{\max_\alpha |\partial_\alpha f(z)|}{\min_\alpha |\partial_\alpha f(z)|}$, and we find that $f$ is $M^2$-QC. This implies that the refined charts are contained in a unique QC-structure on $X$. The last statement is proved analogously.
\end{proof}

\begin{remark}
Note that on a fixed Hausdorff compact topological space, different smooth structures may correspond to the same QC structure. To see this, let $X$ be a compact smooth manifold with $\partial$, let $\varphi : X \to X$ be a self-homeomorphism which is QC but not smooth\footnote{It is trivial to construct homeomorphisms which are $C^1$, hence QC, but not smooth. Probably the easiest way to construct a QC mapping which is not even $C^1$ is by starting from a $C^1$-diffeomorphism $\varphi : S^1 \to S^1$ and creating a QC mapping with a singularity point as in the proof of lemma \ref{lemma:corners}. More elaborate singularities along smooth arcs can be constructed by gluing two smooth maps, cf. lemma \ref{lemma:local pushout}.}. We can define a second, different, smooth structure on $X$ simply by pulling back the smooth structure along the homeomorphism $\varphi$. Both smooth structures correspond to the same QC structure on $X$.
\end{remark}

\section{Beltrami Differentials}
\label{subsec:beltrami-diffs}
In what follows $X$ denotes a fixed QC-surface. We want to give a managable description of all conformal structures on $X$ which are compatible with the QC structure. As we remarked in chapter \ref{subsec:analyticdef}, the complex dilatation is a good way to encode almost complex structures relative to a fixed complex structure. Beltrami differentials can be thought of as an invariant way of doing the same thing, ridding ourselves of the reference complex structure. We now make this more precise.

Choose a finite QC atlas for $X$, $(\varphi_i)_{i =1}^n$. We can then define
\textbf{the bundle of beltrami differentials}, $\mbox{Belt}(X) \to X'$. It is a bundle defined above $X' = X - E$, where $E$ is some subset of measure zero (we will explain precisely what that means below) which depends on the atlas we have chosen.

Let $\psi_{ij} = \varphi_j \circ \varphi_i^{-1}$ denote the QC transition maps. The partial derivatives $(\psi_{ij})_z, (\psi_{ij})_{\overline z}$ exist except on a set of measure zero $E_{ij} \subset V_{ij}$. We can also assume these sets are symmetric: $\psi_{ij}(E_{ij}) = E_{ji}$. Here we are using the fact that QC maps preserve sets of measure zero with respect to the Lebesgue measure on $\overline \hh$ (see \citepalias[pg. 165]{LV}). This fact also implies that we can define a measure zero set on $X$ as a set whose images under all QC charts are measure zero. Letting $E_i := \cup_j E_{ij} \subset V_i$, we find that $E = \cup \varphi_i^{-1}(E_i)$ is such a set.

We define $\belt(X) := \left(\coprod (V_i - E_i) \times B_1\right) / \sim$. Here $B_1 := \{z \in \cc | |z| < 1\}$, and we set
\begin{equation}
\label{eqn:gluing}
(z,\mu)_i \sim \left(\psi_{ij}(z), \frac{\psi_z}{\overline{\psi_z}} \frac{\mu - \mu_\psi}{1 - \overline{\mu_\psi} \mu}\right)_j\end{equation}
for any $z \in V_{ij} - E_i$ (cf. \textsection \ref{subsec:transformation}) We have a natural projection $\belt(X) \to X' := X - E$.

The bundle constructed above depended on the choice of finite atlas. Different choices will lead to bundles which are uniquely isomorphic when appropriately restricted: if $\mathcal{A}_1, \mathcal{A}_2$ are two finite atlases then so is $\mathcal{A}_{12}:=\mathcal{A}_1 \cup \mathcal{A}_2$, and we can take a common $E$ for all three bundles. Then the maps $[(z,\mu)_i)]_{\mathcal{A}_j} \mapsto [(z,\mu)_i]_{\mathcal{A}_{12}}$, $j = 1,2$ are isomorphisms of bundles over $X - E$. The above discussion leads naturally to consideration of measurable sections of $\belt(X)$ defined ``up to measure zero modification''. We denote these by $\Gamma(X)$. This set \emph{is} canonically defined, in the sense that there's a unique isomorphism between the $\Gamma$'s constructed from any two choices of finite atlases.

We can endow $B_1$ with the hyperbolic metric; it is invariant under the action of $PSL_2(\rr)$ given by eq. (\ref{eqn:gluing}), and in this way we obtain a metric (the \textbf{Teichm\"uller metric}) on the fibers given by
\[d(\mu_1,\mu_2) = \log \left(\frac{1 + \delta}{1 - \delta}\right), \delta = \left|\frac{\mu_1 - \mu_2}{1 - \mu_1 \overline{\mu_2}}\right|.\]

We would like to define a metric on $\Gamma(X)$ using the essential supremum of the fiberwise distance. But it is not hard to construct examples of $\mu_1,\mu_2$ where the essential supremum is infinite.

To remedy this we want to consider the \textbf{bounded Beltrami differentials on $X$}, denoted $\Gamma_\infty(X) \subset \Gamma(X)$. The definition again involves a choice of finite atlas, but in a non-essential way. We say that $\mu$ is bounded if $\|\mu_i \|_\infty < 1$ for $1 \leq i \leq n$, or equivalently, if $\|D_{\mu_i}\|_\infty < \infty$ where $D_{\mu_i}(z) = \frac{1 + |\mu_i|}{1 - |\mu_i|}$ is the dilatation of $\mu$. This definition does not depend on the choice of atlas, for dilatations are multiplied under composition and the transitions are finite in number and have bounded dilatation.

It is also not hard to construct a bounded Beltrami differential by induction, call it $\mu$. The bounded differentials are then precisely those which are a finite distance from $\mu$.

For any QC-morphism $f : X \to Y$ we have a map $\Gamma_\infty(f) = f^* : \Gamma_\infty(Y) \to \Gamma_\infty(X)$. If $\mu \in \Gamma_\infty(Y)$ admits a local representation in terms of the measurable function $\nu_Y$ in some local variable $w = \varphi_Y$ on $V \subset Y$ then $f^* \mu$ is given in terms of a local variables $z = \varphi_X$ on $X$ by the function $\nu_X$,
\begin{equation}
\label{eqn:pullback}
\nu_X = \frac{1}{1 - |\mu_g|^2}\frac{(\nu_Y \circ g) + u \mu_g}{u + \overline {\mu_g} (\nu_Y\circ g) },
\end{equation}

where $g = \varphi_Y \circ f \circ \varphi_X^{-1}$ and $u = \frac{g_z}{\overline{g}_{\overline z}}$ (this is just the inverse of the transformation rule).

This turns $\Gamma_\infty$ into a contravariant symmetric monoidal functor (cf. \textsection \ref{subsec:monoidal}) from (QC-surfaces,$\coprod$) to (metric spaces, $\times$), where $d_{X \times Y}((x,y),(x',y')) = \max(d_X(x,x'),d_Y(y,y'))$ .

\section{Conformal Structures on a QC-surface}
\label{subsec:conformal structures}
We would like to show that

\vspace{2mm}
\noindent
\emph{$\Gamma_\infty(X)$ can be canonically identified with the set of conformal structures on $X$ which are QC equivalent to $X$.}
\vspace{2mm}

Fix $\mu \in \Gamma_\infty(X)$, and for every chart $\varphi_i : U_i \to V_i$, $1 \leq i \leq n$ define $\psi_i = \omega_{\mu_i} \circ \varphi_i$. Here $\omega_{\mu_i} =: \omega$ is a solution (in the sense of theorem \ref{thm:integration}) of $\omega_z = \mu_i \omega_{\overline z}$, where $\mu_i$ is the representation of $\mu$ as a section of $B_1 \times V_i$, extended by $0$ to the half plane and then by reflection, $\mu(\overline z) = \overline{\mu(z)}$ to all of $\cc$. We normalize the solution so that it fixes $0,1$ and $\infty$.

Let us show that $\psi_i$ maps $U_i$ to an open subset of $\overline \hh$. If we denote by $C$ reflection about the real axis, we find that $\omega_{\mu_i}$ and $C \circ \omega_{\mu_i} \circ C$ have the same complex dilatation, and hence are equal up to a M\"obius transformation $A$. $A$ must fix $0,1,\infty$, and is therefore the identity. We can now conclude that $\omega$ maps the real line onto the real line, and the upper half plane to itself, and so, given that the image of $\varphi_i$ is contained in the closed upper half plane we conclude that the image of $\psi_i = \omega \circ \varphi_i$ is also contained in it, as desired.

A straightforward computation shows that the transition maps $\psi_j \circ \psi_i^{-1}$ are conformal (indeed, the gluing condition (\ref{eqn:gluing}) was designed for this purpose), and so the $\psi_i$ determine the structure of a Riemann surface (with $\partial$) on $X$.

We still need to show that the map
\[\Gamma_\infty(X) \to \{\mbox{conformal, QC-equivalent structures on $X$}\}\] is bijective; it is useful to phrase it in slightly more general terms.

\begin{definition}
\label{def:pullback}
Assume that $Y$ is given a Riemann surface structure, e.g., a conformal structure defined by some $\mu \in \Gamma_\infty(Y)$. We can define $\mu_f \in \Gamma_\infty(X)$ by setting $\left(\mu_{f}\right)_i(z) = \mu_{\psi_{f(z)} \circ f \circ \varphi_i^{-1}}$, where $\psi_{f(z)}$ is some conformal chart of a neighbourhood of $f(z)$.
\end{definition}

Replacing $\psi$ by a conformally equivalent $\psi'$ does not change $\mu$ so this gives a well-defined local section. We need to verify that these sections glue together according to (\ref{eqn:gluing}). Indeed:
\[(\mu_f)_j \circ \psi_{ij} = \mu_{\psi \circ f \circ \varphi_j^{-1}} \circ \psi_{ij} = \mu_{(\psi \circ f \circ \varphi_i^{-1}) \circ \psi_{ij}^{-1}} \circ \psi_{ij} = (\mu_f)_i\]\qed

We can now define the inverse mapping
\[\{\mbox{conformal, QC-equivalent structures on $X$}\} \to \Gamma_\infty(X)\] simply by sending the conformal structure to $\mu_{id}$, where $id : X \to X$ is the identity map between the QC surface $X$ and the conformal surface which is also $X$.

\vspace{3mm}
Definition \ref{def:pullback} allows us to define a pull back map $\Gamma_\infty(f) : \Gamma_\infty(Y) \to \Gamma_\infty(X)$ for $f : X \to Y$ a QC morphism, by $\mu \mapsto \mu_f$; here we are using $\mu$ to define the complex structure on $Y$, and then appealing to the definition of $\mu_f$ above. This map corresponds to the pull-back of conformal structures from $Y$ to $X$ along $f$; we have given an explicit formula for this pullback in the previous section, see eq. (\ref{eqn:pullback}).

\vspace{3mm}
We remark that if $f : X \to Y, g : X \to Y'$ satisfy $\mu_f = \mu_g$ then $\mu_{\psi' \circ g^{-1} \circ f \circ \psi^{-1}} = 0$ for any pair of charts, which shows that $f$ and $g$ differ by a conformal map $Y \to Y'$.

\section{Open-Closed QC Surfaces}
\label{subsec:OCQC surfaces}
So far we have ignored the parametrization of the boundary. We now take it into consideration.

A \textbf{QC-surface with boundary parametrization between $(C_-,O_-)$ and $(C_+,O_+)$} (or just ``a QC-surface with parametrization'') is a QC-surface with topological embeddings $\alpha_\pm : S_1^{C_\pm} \to \partial X, \beta_\pm : [0,1]^{O_\pm} \to \partial X$ having disjoint images and such that $\alpha_+, \beta_+$ preserve orientation and $\alpha_-, \beta_-$ reverse orientation (cf. definition \ref{def:OCR}).

The parametrization may be smooth with respect to some of the charts, or perhaps with respect to none at all. In accordance with our agreement that OC-surfaces are obtained from OCR-surfaces by forgetting some of the structure, we define an \textbf{Open-Closed-QC surface} to be a QC surface with parametrization such that under \emph{some} conformal atlas compatible with the QC structure the parametrization is smooth. It is natural to seek a more explicit definition of an OC-QC surface.

The following result is central.
\begin{theorem}
\label{thm:QSequiv}
Let $J$ be a compact interval in $\rr$ or $\rr$ itself. Then for a continuous, strictly increasing function $h : J \to \rr$ the following conditions are equivalent:

\begin{enumerate}[(A)]
\item $h$ can be extended to a QC mapping $H : U \to V$ where $U,V$ are neighbourhoods of $\rr \cup \infty$ in $\overline \hh$, and $H_{|J} = h$.
\item there exists some positive $k$ such that
\[1/k \leq \frac{h(x + t) - h(x)}{h(x) - h(x-t)} \leq k\] for all $x$ and $t$ such that $x-t,x,x+t \in J$.
\end{enumerate}
\end{theorem}

\begin{proof}
The theorem is a straightforward combination of results from Part II of \citepalias{LV}:
\begin{itemize}
\item Theorem 6.2 and Theorem 6.3, which characterize functions satisfying (B) in $\rr$ as maps which are the boundary values of QC self-mappings of the open half plane.
\item Lemma 7.1, which says that a function which satisfies (B) on a compact interval can be extended to one which satisfies it in $\rr$.
\item Theorem 6.4, which says that the requirement in (A) is equivalent to the existence of a self-mapping of $\overline \hh$ which restricts to $h$ on $\rr$.
\end{itemize}
\end{proof}

A function is called \textbf{quasisymmetric} if it satisfies any of the conditions above. We now show that the appropriate extension of this notion to boundary parametrizations of QC-surfaces will give the desired characterization of OC-QC surfaces. Roughly speaking, a map $\gamma : D \to \partial X$ ($D = S^1$ or $D = [0,1]$) is quasisymmetric if its image under a collar charts which send the boundary circle to $\rr \cup \infty$ is quasisymmetric. We now formulate this more precisely.

Let $X$ be a QC surface with boundary parametrization, and focus on one of the boundary circles $C$. \textbf{A collar chart for $C$} is a surjective QC mapping $\varphi : U \to A$ such that:

\begin{enumerate}
\item $U$ an open neighbourhood of $C$,
\item $A$ is the extended closed upper half-plane minus a disk, and
\item $\varphi$ takes $C$ to $\rr \cup \infty$.
\end{enumerate}

To construct a collar chart for $C$ fix some auxiliary complex structure on $X$ compatible with the QC structure. It is a well-known result of Morse theory that there exists an open neighbourhood of the boundary circle which is diffeomorphic to $S^1 \times [0,1]$. Let $U$ be one such neighbourhood. By the uniformization theorem $U$ is conformally equivalent to an annulus in the complex plane. A suitable M\"obius transformation then maps this annulus to $A$ with the boundary circle mapped to the extended real line. Note that by rotating $\varphi$ we may choose freely which point on the boundary circle is mapped to $\infty$.

Let $\gamma : D \to \partial X$ be a component of the parametrization, where $D = [0,1]$ for open and $D = \rr \cup \{\infty\} \simeq S^1$ for closed boundary components. Let $\varphi$ be any collar chart for the target boundary circle such that the preimage of $\infty$ is $\gamma(\infty)$ if $D = \rr \cup \{\infty\}$ or some point outside $\gamma(D)$ if $D = [0,1]$. If $\gamma$ is a component of the orientation-preserving parametrizations $\alpha_+,\beta_+$ we say it is quasisymmetric if $\varphi \circ \gamma_{|\rr} : D \cap \rr \to \rr$ is quasisymmetric. For components of $\alpha_-, \beta_-$ we require that $\varphi \circ \gamma \circ \tau_{|\rr}$ be quasisymmetric, where $\tau : x \mapsto -x$. This does not depend on the choice of collar chart since the transition between different charts is a QC mapping which can be extended to a neighbourhood of $\rr \cup \infty$, and condition (A) is clearly invariant under such transitions. The boundary parameterization is called quasisymmetric if every component of the parametrization is quasisymmetric.

We can now prove
\begin{corollary}
The following conditions are equivalent:

\begin{enumerate}[(a)]
\item $X$ is an OC-QC surface.
\item $X$ is a QC-surface with quasisymmetric boundary parametrization.
\end{enumerate}
\end{corollary}

\begin{proof}
(a) $\Rightarrow$ (b): Let $\gamma$ be a component of the parametrization. We assume that $X$ admits an OC-Riemann surface structure; this means that if we construct a collar chart $\varphi$ which is conformal with respect to this structure the map $h := \varphi \circ \gamma_{|\rr} : D \cap \rr \to \rr$ will be a smooth embedding. For $D = [0,1]$ we find that $h'$ obtains a maximum $M$ and a minimum $m$ in $D$, both positive, and then, by the mean value theorem, $\frac{m}{M} \leq \frac{h(x + t) - h(x)}{h(x) - h(x-t)} \leq \frac{M}{m}$, which shows that $h$, and hence $\gamma$, is quasisymmetric. If $D = \rr \cup \infty$ it is still true that $h'$ is bounded, since it is a diffeomorphism also ``at $\infty$''\footnote{More precisely, if $h$ is a diffeomorphism at $\infty$ then we have $h(x) = \frac{1}{c\frac{1}{x} + \alpha(1/x)}$ where $c > 0$ is the derivative at $x = \infty$ relative to the parameter $y = 1/x$, and $\alpha$ is the remainder term - a continuously differentiable function with $\frac{\alpha(y)}{y} = x\alpha(1/x) \to 0$ as $y \to 0$ (or $x \to \infty$). Differentiating we find that \[h'(x) = \frac{1}{c + x\alpha(1/x)} - \frac{x\alpha(1/x)}{\left(c + x\alpha(1/x)\right)^2} - \frac{x^2\alpha'(1/x)\frac{-1}{x^2}}{\left(c + x\alpha(1/x)\right)^2}.\] So as $x \to \infty$, the RHS tends to $\frac{1}{c}$; in particular the derivatives are bounded.}.

(b) $\Rightarrow$ (a):
Choose any Beltrami differential $\mu_0$ on $X$. We want to show that it can be modified near the boundary so that the parmeterization becomes smooth. Let $\gamma_1 : D \to \partial X$ be one component of the parametrization (as before, $D = S^1$ or $D = [0,1]$). Choose some QC chart $\varphi$ of a collar neighbourhood of the boundary circle containing the image of $\gamma$. $\varphi \circ \gamma$ is quasisymmetric, and so by theorem \ref{thm:QSequiv} we can extend it to a QC mapping $\psi$ of a neighbourhood of $\rr \cup \infty$. Let $\widetilde \varphi = \psi^{-1} \circ \varphi$; it is a quasiconformal chart of a (perhaps smaller) collar neighbourhood, such that $\widetilde \varphi \circ \gamma$ is the standard embedding of $D$ into $\cc$ (in particular, it is smooth). We now use $\widetilde \varphi$ to \emph{define} a complex structure near the image of $\gamma$. Explicitly, let $U_1 \subset X$ be a small open set containing the image of $\gamma$, and contained in the domain of $\widetilde \varphi$. We can then define $\mu_1$ to be the same as $\mu_0$ outside $U_1$ and the section represented by 0 in the chart $\widetilde \varphi$ in $U_1$. This ensures that the image of $\gamma_1$ is smooth (in fact, analytic) relative to the conformal structure $\mu_1$.

\begin{figure}[hb]
  \centering
\includegraphics[width = 100mm]{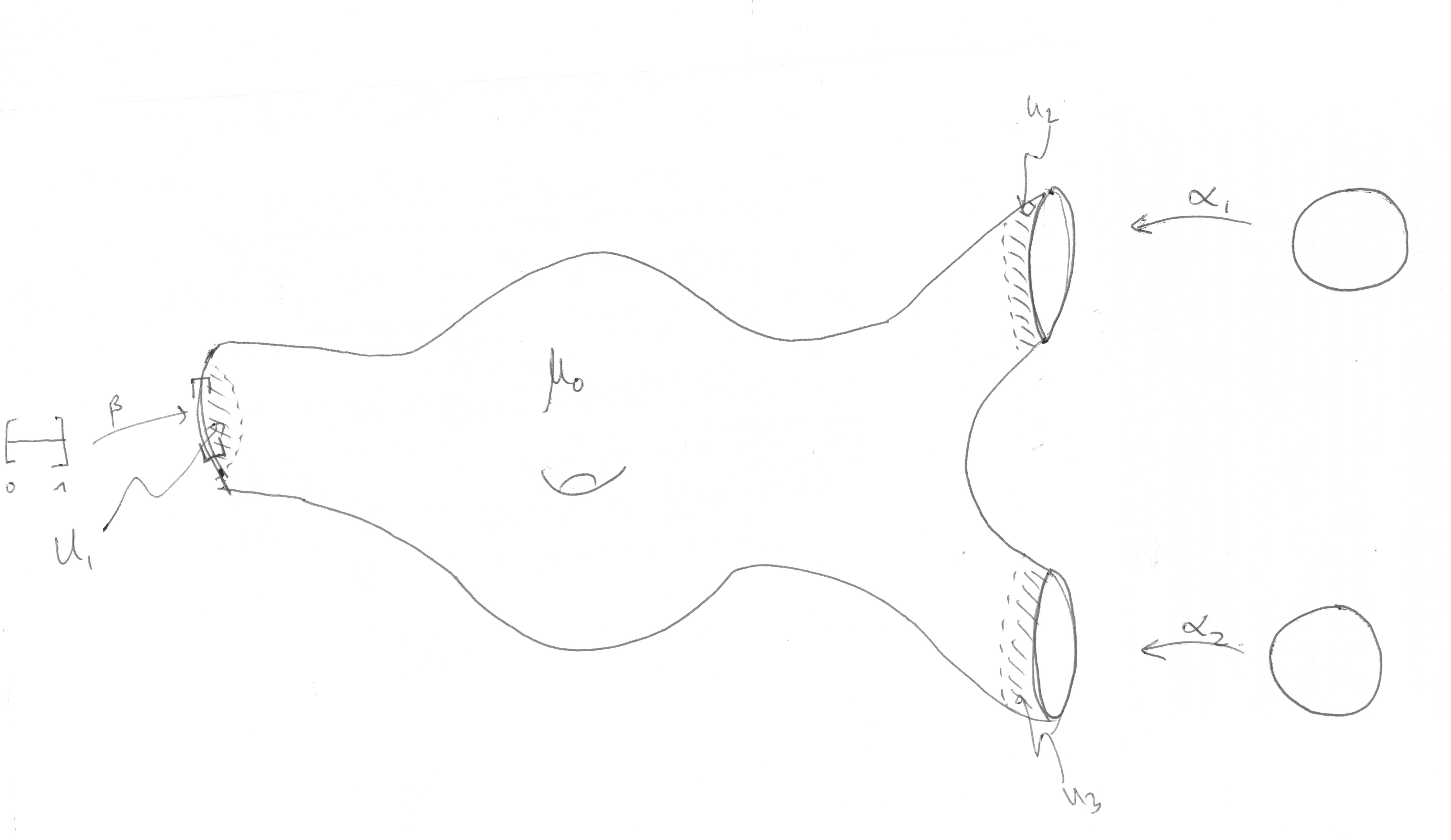}
  \caption[modifying]
  {Modifying the Beltrami differential near the parametrized boundary.}
\end{figure}

If $\mu_n$ has been constructed so that the parametrization components $\gamma_1,...,\gamma_n$ are smooth near the boundary relative to $\mu_n$, then we can follow the same procedure to construct $\mu_{n+1}$; we just need to be careful to choose $U_{n+1}$ small enough so that its closure does not contain any points in the image of $\gamma_1,...,\gamma_n$.
\end{proof}

\begin{remark}
For an OC-QC surface, it is natural to ask which Beltrami-differentials $\mu$ represent a complex structure relative to which the parametrization is smooth. This is clearly determined by the ``germ at the parametrized boundary'' (we consider sections up to an equivalence which identifies any two sections which agree on an open neighbourhood of the image of the parametrization). If we fix a smooth structure on $X$ relative to which the parametrization is smooth, then a sufficient condition is that $\mu$ is smooth near the boundary in every smooth chart. This follows from the fact that the Beurling-Ahlfors extension (see \citepalias[II \textsection 6.5 ]{LV}) is smooth for smooth $\mu$.

This condition is not necessary. To see this, note that we can modify the smooth structure on $\overline \hh$ so that the imaginary axis is a ``fault line'' - we identify $Q_1$ and $Q_2$ so that the constant vector field $\hat x$ on $Q_1$ is identified with $(1 + e^{-1/y^2})\hat x$ on $Q_2$. The set of smooth ``parametrizations'' $\rr \to \rr$ is not modified by this change, but a function $\overline \hh \to \rr$ which is smooth with respect to one structure is not smooth with respect to the other, even if we restrict to an arbitrarily small neighbourhood of the real axis.

It should be interesting to try and refine this in order to obtain a necessary and sufficient condition on the germ of the Beltrami differential.
\end{remark}

\chapter{Morphism Spaces and Composition Maps}
\label{chapt6}
\label{sec:composition}
We are now ready to tackle the problem of defining continuous composition maps.
We first redefine the morphism space in the language of QC-surfaces.

\section{The Morphism Space}
\label{subsec:the morphism space}
We have seen that an OC-smooth structure induces the structure of a QC-surface with quasisymmetric parametrization; that a diffeomorphism rel-$\partial$ induces a QC-isomorphism rel-$\partial$; and that if the parametrization is quasisymmetric then there exists a compatible OCR-surface (which is in particular an OC-smooth surface). This can be summarized by saying that the forgetful functor
\[\mbox{OC-smooth} \to \mbox{OC-QC}\]
between these two categories of surfaces and isomorphisms rel-$\partial$ between them is faithful and essentially surjective. It is also ``essentially injective'', that is, 1-1 on isomorphism types: In appendix \ref{sec:OCsmoothtype} we show that OC-smooth (isomorphism) types are determined by topological data. QC-isomorphisms rel-$\partial$ clearly preserve this kind of data, so the forgetful functor is 1-1 on types.

Note, however, that this functor is \emph{not} an equivalence of categories, since it is not full: there are ``more'' QC-morphisms. In fact it is precisely this added flexibility that we are looking for.

\vspace{2mm}
In chapter \ref{sec:smoothdef}, $\Sigma_\tau$ represented a fixed OC-smooth model surface for the OC-smooth type $\tau$. From now on we will think of $\Sigma_\tau$ merely as an OC-QC surface. That is, a QC-surface with fixed quasisymmetric parameterizations.

We write $\widetilde \mm_\Lambda(\tau) := \Gamma_\infty(\Sigma_\tau) / QC(\Sigma_\tau)$ for \textbf{the quasisymmetric morphism space of type $\tau$}. Here $QC(\Sigma_\tau)$ are all QC automorphisms rel-$\partial$ of $\Sigma_\tau$ (i.e., all QC self-isomorphisms which fix the parameterized part of the boundary). These act on $\Gamma_\infty(\Sigma_\tau)$ as explained in \textsection \ref{subsec:conformal structures}.

The same arguments as in chapter \ref{sec:smoothdef} show that this set is in bijection with the set of ``Quasisymmetric Riemann Surfaces'' of type $\tau$ modulo equivalence. These are the same as OCR-surfaces except that the parametrization need not be smooth but merely quasisymmetric. While this seems a very natural space to consider, \citep{segal}, \citep{costello} focus on \textbf{the smooth morphism space of type $\tau$}, $\mm_\Lambda(\tau) \subset \widetilde \mm_\Lambda(\tau)$, which consists of those (equivalence classes of) conformal structures relative to which the parameterization is smooth.

The reuse of the notation $\mm_\Lambda(\tau)$ from chapter \ref{sec:smoothdef} is justified since both spaces are in a natural bijection, as one verifies readily\footnote{In fact, we strongly believe that this map is a homeomorphism, although continuity in one of the directions requires some subtle arguments with Teichm\"uller spaces which need to be ironed out.}.

The total morphism space is defined by the topological disjoint union over types:

\[\mm_\Lambda((C_-,O_-),(C_+,O_+)) := \coprod_{\{\tau | \tau : (C_-,O_-) \Rightarrow (C_+,O_+)\}} \mm_\Lambda(\tau) \]

\section{Gluing OC-QC Surfaces}
We now turn to the problem of defining continuous composition maps
\[\mm_\Lambda((C_-,O_-),(C_0,O_0)) \times \mm_\Lambda((C_0,O_0),(C_+,O_+)) \to \mm_\Lambda((C_-,O_-),(C_+,O_+))\]

The composition of morphisms is obtained by the sewing of representative OCR-surfaces.

In \textsection\ref{subsec:sewing sketch} we have sketched how to define this sewing for any pair of composable surfaces. However, if we want to show that this process induces a continuous map of the morphism spaces we must be more organized, and define a sewing map by first gluing the underlying QC-surfaces. In the following we use results from \textsection \ref{subsec:bending and gluing} to prove the existence of such a gluing and prove it is unique up to unique isomorphism.

\begin{lemma}
\label{lemma:pushout}
Let $\tau_1$ and $\tau_2$ be composable types of OC-QC surfaces, with corresponding models $\Sigma_{1}$ and $\Sigma_{2}$. There exists a structure of an OC-QC surface on the topological pushout $\Sigma_{12}$, so that $\iota_1,\iota_2$ become QC embeddings.

\hspace{30mm}\xymatrix{
(S^1)^C \times [0,1]^O  \ar[d]^{\alpha^2_- \times \beta^2_-} \ar[r]_{\alpha^1_+ \times \beta^1_+}
& \Sigma_{1} \ar[d]_{\iota_1} \ar@/^/[ddr] & \\
\Sigma_{2} \ar@/_/[rrd] \ar[r]_{\iota_2} & \Sigma_{12} \ar@{.>}[dr]|-{\exists !}& \\
& & Y
}

\noindent
This QC structure satisfies the universal property of the pushout: for any pair of morphisms $\Sigma_1 \to Y$ and $\Sigma_2 \to Y$ making the outer parametrization square commute there exists a unique QC morphism $\Sigma_{12} \to Y$ making the side triangles commute.
\end{lemma}

\begin{proof}
We construct a chart for every point $p \in \Sigma_{12}$.
We use the following terminology: the common parametrization refers to the outgoing parametrization components of $\Sigma_1$ and the incoming parametrization components of $\Sigma_2$. We call the image of the common parametrization in $\Sigma_{\tau_j}$ ``the common boundary'' (so the common boundary consists of the points which are glued). The image of the common boundary in $\Sigma$ is called ``the seam''. We now distinguish between three cases:

\vspace{2mm}
\noindent
\emph{Case I: $p$ is not on the seam.} That is, $\iota_j^{-1}(p)$ is not in the image of the common parametrization. In this case it belongs to exactly one piece, say $\Sigma_j$, and we take a smooth chart of $\Sigma_j$ which avoids the seam.

\vspace{2mm}
\noindent
\emph{Case II: $p$ is on the seam, but not on the boundary of $\Sigma$.} This means that $p_j := \iota_j^{-1}(p)$ belongs to the common boundary for $j = 1,2$, but is not the image of an open boundary endpoint in either piece. Let's assume that $p_1 = \beta^1_+(t)$ and $p_2 = \beta^2_-(t)$ for some $t \in (0,1)$ (the proof for a closed boundary component is similar). There exist two smooth charts, $\varphi_j : \Sigma_j \supset U_j \to V_j$, where $U_j$ is an open neighbourhood of $p_j$, and $V_1,V_2$ are open subsets of the upper and lower half-planes, resp..

Choose some $\delta > 0$ so small that $(t-\delta, t+\delta) \subset [0,1]$ and so that the image of this interval under $\beta^{j}_\pm$ is contained in $U_j$. Let $I_j = \varphi_j(\beta^{j}_\pm((t-\delta,t+\delta)))$. Restricting as necessary we may assume $V_j \cap \rr = I_j$ and that $V_j$ is contained in the vertical strip lying above (resp., below) $I_j$.

We may apply the diffeomorphism $\rho : (x,y) \mapsto (\beta^1_+((\beta^2_-)^{-1}(x)),y)$ to $V_2$, thus aligning it with $V_1$. The desired chart is then $\varphi_1 \cup \rho \varphi_2 : U_1 \cup U_2 \to V_1 \cup \rho(V_2)$.

\begin{figure}[hb]
  \centering
\includegraphics[width = 60mm]{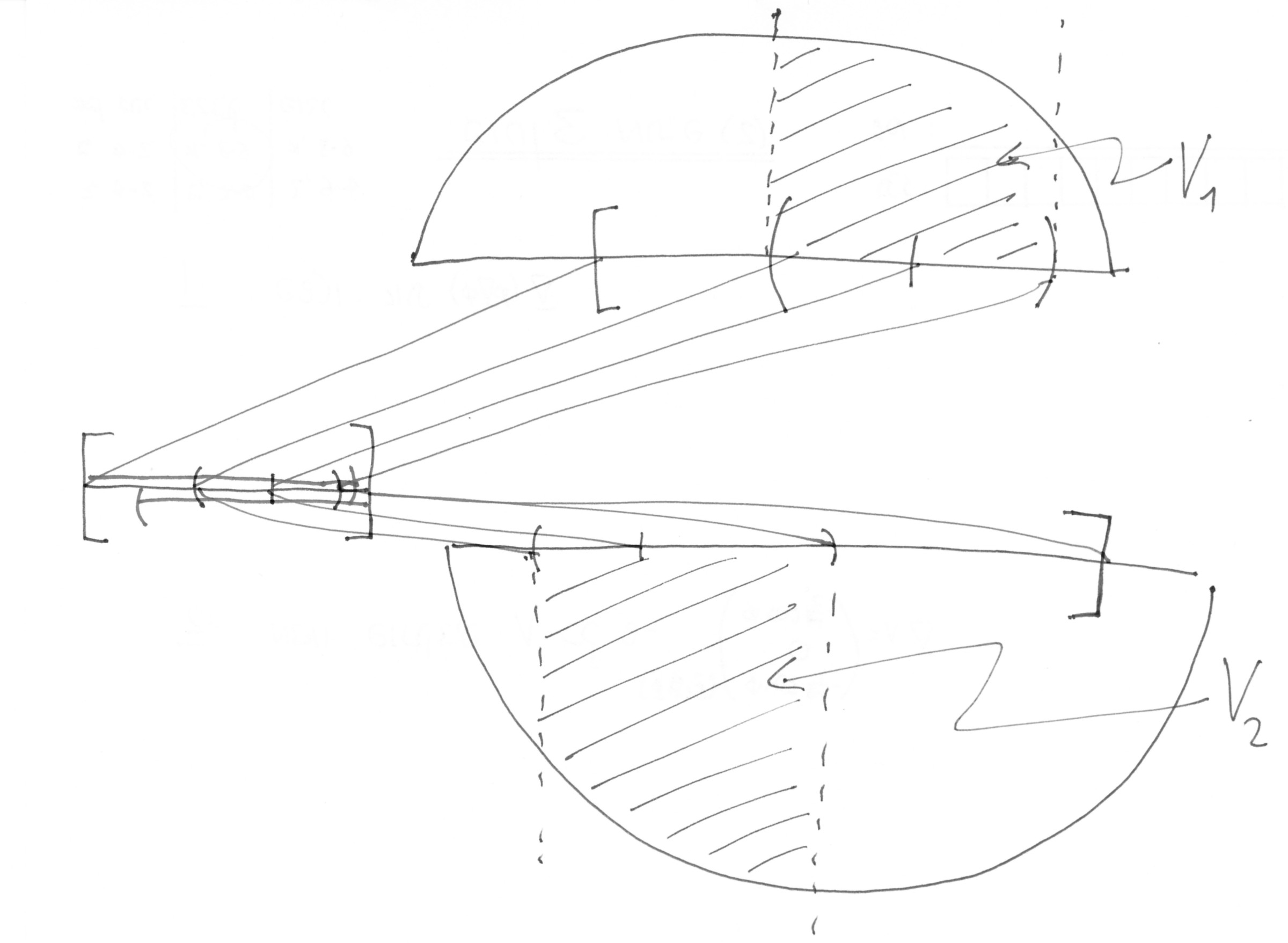}
\includegraphics[width = 60mm]{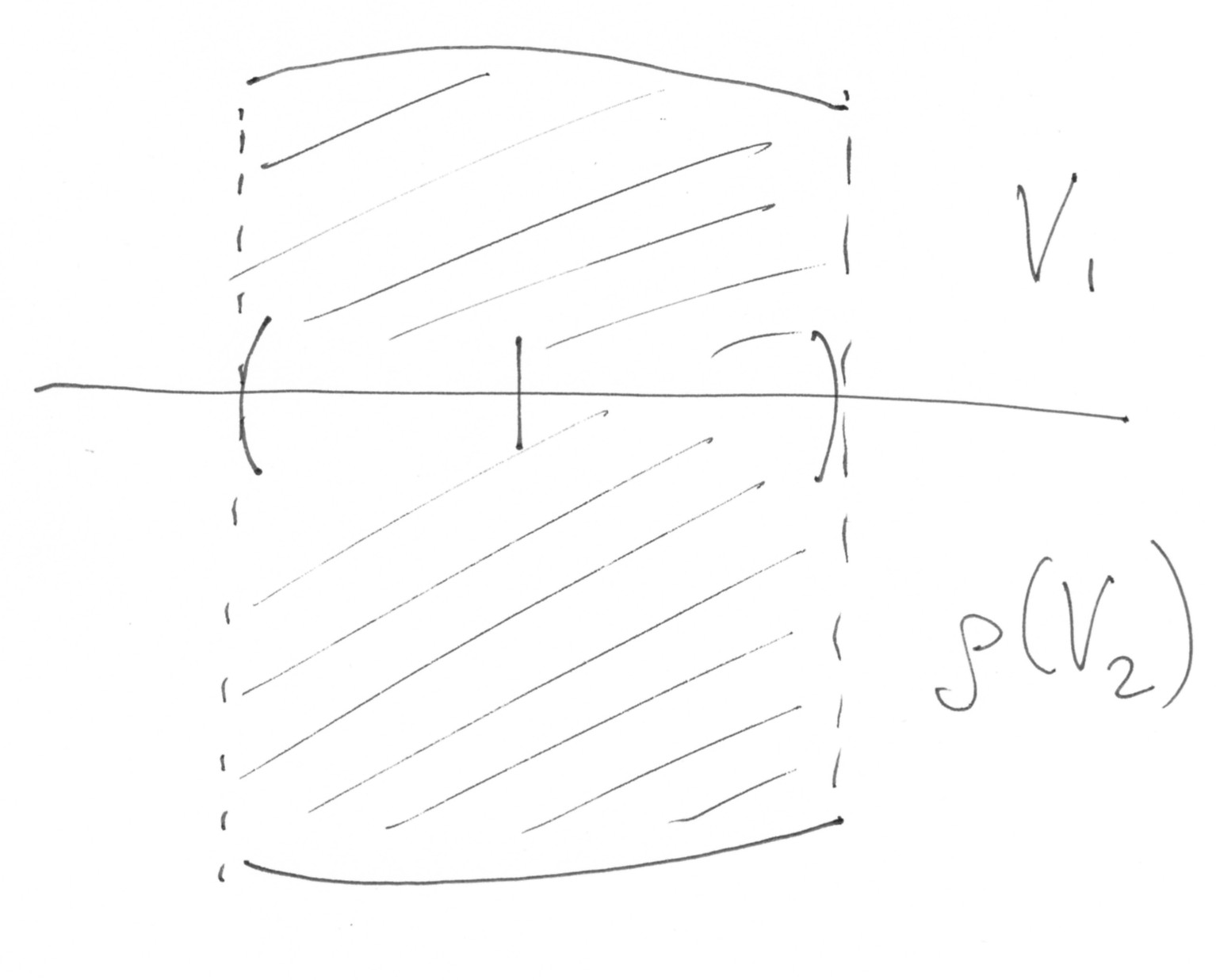}
  \caption[aligning]%
  {A smooth map is used to align the two charts on the left according to the reparametrization, in a neighbourhood of an interior point of the open boundary segment. The resultant two charts on the right can then be juxtaposed to obtain a chart of the sewed surface.}
\end{figure}

\vspace{2mm}
\noindent
\emph{Case III: $p$ is the image of an endpoint of a common open boundary segment.}
We may assume $p_1 = \beta^1_+(0), p_2 = \beta^2_-(0)$. Using a similar argument to the one used before (restricting and aligning using a ``horizontal'' diffeomorphism) we may assume that the domain of both charts intersects the parameterized boundary only in $\beta^j_\pm([0,\delta))$ for some $\delta > 0$, and that $\varphi_1(\beta^1_+(t)) = \varphi_2(\beta_-(t))$ for every $t\in [0,\delta)$.

For convenience we set $\varphi_1(\beta^1_+(0)) = \varphi_2(\beta^2_-(0)) = 0$. Note that the points whose images lie in the negative real axes are not to be glued. We want to split open these two rays so they become a smooth boundary for the glued surface. To achieve this we apply lemma \ref{lemma:corners} and construct $\sigma_1, \sigma_2$ which map the upper (resp., lower) closed half plane to the first (resp. second) quadrant.

As a chart for $p$, take $\sigma_1 \varphi_1 \cup \sigma_2 \varphi_2$.

This completes the construction of the atlas. One verifies that the transition maps are indeed QC. By construction, $\iota_1,\iota_2$ become QC embeddings.

We note that if the parametrizations of the pieces were quasisymmetric then so is the resultant parametrization of the whole, which shows that $\Sigma_{12}$ is an OC-QC surface.

Finally, we need to check the universal property holds. Let $\varphi_i : \Sigma_i \to Y, i=1,2$ be a pair of morphisms which make the outer square of the diagram commute. By the universal property of the topological push-out, there exists a unique continuous map $\varphi = \varphi_1 \cup_{\alpha,\beta} \varphi_2$ which makes the triangles commute; in-fact this map must be a homeomorphism or an embedding. We may verify that it is QC by working with the atlas we've just constructed and applying lemma \ref{lemma:local pushout}.
\end{proof}

\begin{definition}
An atlas as in the proof above will be called \textbf{aligned}. More precisely, if $\Sigma_{12}$ is the push-out of $\Sigma_1,\Sigma_2$, an aligned atlas is one in which points of the seam are mapped to the real axis or the imaginary axis, and an endpoint of an open boundary is mapped to the origin.
\end{definition}

%We are now ready to define the composition. We cannot avoid some abstract category-theoretical discussion, but to avoid obscuring the essentials we choose to postpone it to the next chapter. We assume now that $\Sigma_{123}$ is some smooth surface, and that we cut into three pieces along smooth curves whose endpoints, if any, lie on the boundary, so that $\Sigma_{123}$ decomposes into three OC-smooth surfaces $\Sigma_{123} = \Sigma_1 \cup \Sigma_2 \cup \Sigma_3$ which QC-embed into it.
%
%\begin{theorem}
%There is a continuous map $\Gamma_\infty(\Sigma_1) \times \Gamma_0(\Sigma_2) \to \Gamma_0(\Sigma_1 \cup \Sigma_2)$
%
%
%$\Gamma_0(\Sigma_2) \times \Gamma_0(\Sigma_3) \to \Gamma_0(\Sigma_2 \cup \Sigma_3)$, $\Gamma_0(\Sigma_1 \cup \Sigma_2) \times \Gamma_0(\Sigma_3) \to \Gamma_0(\Sigma_{123})$ and $\Gamma_0(\Sigma_1) \times \Gamma_0(\Sigma_2 \cup \Sigma_3) \to \Gamma_0(\Sigma_{123})$ which are associative in the obvious sense.
%\end{theorem}
%
%\begin{proof}
%We assume the atlas
%\end{proof}
%

\section{Sewing Complex Structures ...}
We're now ready to give a rigorous definition of a continuous sewing operation.
%Let $\iota_1 \coprod \iota_2 : \Sigma_1 \coprod \Sigma_2 \to \Sigma_{12}$ be the

\begin{lemma}
$\Gamma_\infty(\iota_1 \coprod \iota_2) : \Gamma_\infty(\Sigma_{12}) \to \Gamma_\infty(\Sigma_1) \times \Gamma_\infty(\Sigma_2)$ is invertible.
\end{lemma}

\begin{proof}
Clearly $\Gamma_\infty(\iota_1 \coprod \iota_2)$ is an isometric embedding, and in particular 1-1. We need to show that it is surjective.

Let $\mu_j \in \Gamma_\infty(\Sigma_j)$. Assume we have an aligned atlas for $\Sigma_{12}$, indexed by $\alpha$, and define $(\mu_{12})_\alpha$ by the values of $\mu_1$ and $\mu_2$ in the subordinate charts of $\Sigma_1$ and $\Sigma_2$ in the obvious way. Note that near the seam, this involves juxtaposing the sections $(\mu_1)_{\alpha_1}$ and $(\mu_2)_{\alpha_2}$, where $\alpha_j$ is the chart of $\Sigma_j$ obtained by restricting the $\alpha$-chart of $\Sigma_{12}$. These sections  may not agree on the seam, but since it has measure zero\footnote{This is the only place where we are using the fact that the atlas is aligned. In fact, it can be shown that under any QC chart the seam is $\mu$-H\"older continuous for some $\mu < 1$, which implies it has Hausdorff dimension less than 1, and in particular it has measure 0. So we can use any atlas to define $\mu_{12}$, simply juxtaposing the values of $\mu_1$ and $\mu_2$ relative to the restricted charts.} the ambiguity disappears when we consider the result in $\Gamma_\infty(\Sigma_{12})$.
\end{proof}

We denote $\overline{\mathcal{S}} := \Gamma_\infty(\iota_1 \coprod \iota_2)^{-1} : \Gamma_\infty(\Sigma_1) \times \Gamma_\infty(\Sigma_2) \leftrightarrow \Gamma(\Sigma_{12})$
\begin{remark}
\label{rem:sewing complex structures}
In \textsection \ref{subsec:conformal structures} we have seen that $\Gamma_\infty(X)$ correspond to conformal structures compatible with the QC structure on $X$. Suppose $\mu_{12} = \overline{\mathcal{S}}(\mu_1,\mu_2)$. It is straightforward to show that the analytic maps relative to the conformal structure on defined by $\mu_{12}$ are precisely the continuous maps on $\Sigma_{12}$ whose pullbacks along $\iota_1,\iota_2$ are analytic on $\Sigma_1,\Sigma_2$ away from the common boundary, according to the conformal structures defined by $\mu_1$ and $\mu_2$. In appendix \ref{sec:smooth sewing} we show that if the boundary parametrization is smooth with respect to the conformal structures $\mu_1,\mu_2$ (which is always the case for representatives of $\mm_\Lambda(\Sigma_i)$) then the restrictions are in fact smooth up to the boundary.
\end{remark}

\section{... Induces Composition of Morphisms}
Let $\tau_1 : A \to B$ and $\tau_2 : B \to C$ be composable types of OC-surfaces. Let $\Sigma_{1}, \Sigma_{2}, \Sigma_{12}$ be OC-QC models of types $\tau_1, \tau_2$ and $\tau_2 \circ \tau_1$, respectively (cf. appendix \ref{sec:OCsmoothtype}) and let $\iota_j : \Sigma_{j} \to \Sigma_{12}$ be QC-embeddings\footnote{to construct such embeddings it is actually easier to first choose a smooth structure on $\Sigma_{12}$ and cut it along appropriate smooth paths, so that the pieces are diffeomorphic to $\Sigma_1,\Sigma_2$, as we've done in chapter \ref{sec:smoothdef}.}.

\begin{lemma}
The isometry \[\overline {\mathcal{S}} := \Gamma_\infty(\iota_1 \coprod \iota_2)^{-1} : \Gamma_\infty(\Sigma_{1}) \times \Gamma_\infty(\Sigma_{2}) \to \Gamma_\infty(\Sigma_{12})\] descends to a continuous map of morphism spaces:
\[\mathcal{S} : \mm_\Lambda(A,B) \times \mm_\Lambda(B,C) \to \mm_\Lambda(A,C)\]
\end{lemma}

\begin{proof}
By lemma \ref{lemma:pushout} there is a group monomorphism
\[QC(\Sigma_1) \times QC(\Sigma_2) \to QC(\Sigma_{12})\]
(note that the parametrization diagram of lemma \ref{lemma:pushout} commutes since maps in $QC(\Sigma_j)$ are rel-$\partial$). Let $f_j \in QC(\Sigma_j)$, and denote by $f_{12}$ the image of $(f_1, f_2)$ under the group monomorphism. Let $\mu_j \in \Gamma_\infty(\Sigma_j)$. Clearly
\[\Gamma_\infty(f_{12})\overline{\mathcal{S}}(\mu_1,\mu_2) = \overline{\mathcal{S}}(\Gamma_\infty(f_1)\mu_1, \Gamma_\infty(f_2)\mu_2),\]
which shows that there's a well-defined map on the quasisymmetric morphism spaces
\[\widetilde{\mathcal{S}} : \widetilde{\mathcal{M}}_\Lambda(A,B) \times \widetilde{\mathcal{M}}(B,C) \to \widetilde{\mathcal{M}}(A,C)\]

The following diagram shows that $\widetilde{\mathcal{S}}$ is continuous.

\vspace{3mm}
\hspace{30mm}\xymatrix{
\Gamma_\infty \times \Gamma_\infty \ar[d]_q \ar@{.>}[rd] \ar[r]^{\overline{\mathcal{S}}}
& \Gamma_\infty \ar[d]^q \\
\widetilde\mm \times \widetilde\mm \ar[r]^{\widetilde {\mathcal{S}}} & \widetilde \mm }

\noindent
(the dashed map is continuous, hence the bottom map is continuous by the universal property of the quotient topology)

\vspace{3mm}
Finally we restrict to the smooth morphism space (endowed with the subspace topology) obtaining a continuous map
\[\mathcal{S} : \mm_\Lambda \times \mm_\Lambda \to \widetilde{\mm_\Lambda}\]

We claim the image actually lies in $\mm_\Lambda$. This follows from the locality of the smoothness requirement - for every point on the unglued boundary we can take a chart which avoids the seam, and hence is essentially unchanged by the sewing.
\end{proof}

\section{Associativity}
We complete this chapter by proving the associativity of the composition maps.

Note that in order to define the composition maps $\mathcal{S}$ we had to choose, for every composable pair of OC-smooth types $\tau_1,\tau_2$, embeddings $\iota_j : \Sigma_{\tau_j} \to \Sigma_{\tau_2 \circ \tau_1}$. There is no way to make this choice so that the isometries $\overline S$ define an associative operation on the level of sections in $\Gamma_\infty(\Sigma_{\tau_j})$. To see this, take $\tau_1 = \tau_2 = \tau_3$ to be the cylinder with one incoming and one outgoing closed boundary. Clearly $\tau_2 \circ \tau_1 = \tau_3 \circ \tau_2 = \tau_3 \circ \tau_2 \circ \tau_1$ also have the same type, and it is easy to see that no pair of embeddings will make the corresponding maps associative ``on the nose''\footnote{We will replace the morphism space of the cylinder with another topological space shortly; but essentially the same problem appears when composing more complicated surfaces, whose corresponding morphism spaces are not modified further.}.

It turns out, however, that for arbitrary embeddings the induced maps $\overline S$ are associative \emph{up to a (unique) QC isomorphism}, as we now show. This means that when we descend to the morphism spaces the composition maps are associative.

Let $\Sigma_1,\Sigma_2,\Sigma_3$ represent OC-QC models for morphisms \[A \overset{\Sigma_1}\Longrightarrow B \overset{\Sigma_2}\Longrightarrow C \overset{\Sigma_3}\Longrightarrow D.\] The proof is summed up in the following diagram. Straight arrows are QC morphisms, and wiggly ones represent parametrization maps. We construct two morphisms (represented by dashed arrows) using the universal property of the QC pushout.

\vspace{3mm}
\hspace{10mm}\def\g#1{\save
[].[rrrrrr]!C="g#1"*[F.]\frm{}\restore}%
\xymatrix{
\g1 A \ar@{~>}[dr]&& B\ar@{~>}[dl]\ar@{~>}[dr] && C\ar@{~>}[dl]\ar@{~>}[dr] && D\ar@{~>}[dl] \\
&\Sigma_1\ar[dr]&& \Sigma_2\ar[dr]\ar[dl]&& \Sigma_3\ar[dl]&\\
&&\Sigma_{12}\ar[d]\ar@{.>}[drr]^{\exists !} && \Sigma_{23}\ar[d] &&\\
&&\Sigma_{123}\ar@{.>}[rr]^{\exists !} && \Sigma_{123} &&
}

The bottom arrow is invertible, as can be seen by reflecting the diagram and appealing to the uniqueness.

\chapter{Finishing Touches}
\label{chapt7}
\label{sec:finishing}

In this chapter we complete the definition of $\mm_\Lambda$ as a symmetric monoidal category enriched over topological spaces.

These final adjustments follow \citep{costello} closely; in fact, the changes in \textsection \ref{subsec:stability} and \textsection \ref{subsec:D-branes} are somewhat technical and cannot be properly motivated without a better understanding of the theory of moduli and Costello's results. The upside of this is that these details may be skipped at first reading, until a more sophisticated investigation is necessary.

\section{Identity Morphisms}
As defined, the category does not contain identity morphisms. To remedy this, we \emph{replace} the moduli of the cylinder with one incoming and one outgoing closed boundary by $\mbox{Diff}_+(S^1)$ and the moduli of the disc with one incoming and one outgoing open boundary by just a point. More precisely, note that the component of $\mm_\Lambda(A,B)$ which corresponds to a disjoint union of types, $A \coprod A' \overset{\tau \coprod \tau'}\Longrightarrow B \coprod B'$, is equal to the product of components corresponding to $A \overset{\tau}\Longrightarrow B$ and $A' \overset{\tau'}\Longrightarrow B'$, and we have
\[\mm_\Lambda(A,B) = \coprod_{\{\mbox{types $\tau$}\}} \prod_{\{\sigma \subset \tau | \sigma \mbox{ connected type}\}}\mm_{\Lambda,\sigma}\]

So we replace each occurrence of $\mm_\Lambda(\sigma)$, for $\sigma$ a cylinder or a disc as above with $\mbox{Diff}_+(S^1)$ or a point, respectively. There is an obvious way to modify the definition of the composition rule ($\mbox{Diff}_+(S^1)$ acts by changing the parametrization). It is easy to see that associativity is preserved; it is slightly less obvious (but true) that the composition map is still continuous, if we take $\diff(S^1)$ with the $C^1$ topology.

\begin{proof}[Proof Sketch.]
Let $\Sigma_\tau$ be the OC-QC model used to construct $\mm_\Lambda(\tau)$, where $\tau$ has a closed (WLOG, incoming) boundary circle $C \subset \partial \Sigma_\tau$, parametrized with $\alpha_0 : S^1 \to C$. Fix a collar QC chart $\varphi : U \to A$; here $U$ is a collar neighbourhood of $C$ and $A$ is the standard annulus $1 \leq |z| < 2$.

We compose with $\nu \in \diff(S^1)$. Let $\Sigma_\tau'$ be the same as $\Sigma_\tau$, but with $C$ parametrized by $\alpha_0 \circ \nu$. We need to find a QC mapping rel-$\partial$ $\Phi_\nu : \Sigma_\tau \to \Sigma_\tau'$; it is easy to construct one by a self-diffeomorphism $A \to A$ which vanishes near $|z| = 2$ as in the proof of lemma \ref{lemma:smooth twist}. The construction there yields a continuous function $\diff(S^1) \to QC(A \to A)$, where $\diff(S^1)$ is topologized by the $C^1$ norm and $QC(A \to A)$ is topologized by the norm defined by maximal dilatation.

From this it follows that the map
\[\mm_\Lambda(\tau) \times \diff(S^1) \to \mm_\Lambda(\tau)\]
is continuous.
\end{proof}

\begin{remark}
Observe that this map is not uniformly continuous. To wit, note that the module of the collar neighbourhood $U$ depends on the conformal structure $\mu$; the larger the module, the more ``room'' we have to apply $\Phi_\nu$, and the smaller the change in the dilatation (this intuition can be made more precise using the composition formula).
\end{remark}

\section{Stability}
\label{subsec:stability}
We make some additional adjustments to the morphism spaces which exclude Riemann surfaces with an infinite automorphism group.

\begin{enumerate}
\item We require that every connected component of a morphism contains either an incoming closed boundary or a free boundary. This means we simply \emph{erase} those moduli components which do not meet this requirement. Note the remaining components are closed under composition.
\item Discs and annuli which have no parameterized boundaries (only free boundaries) have an infinite group of conformal symmetries. To avoid this we replace the moduli of these two types by a point, in the same way we did for a disc with one incoming and one outgoing open boundary\footnote{It is not clear whether we shall need to record the existence of these components, or just identify morphisms which differ by adding disjoint components of such types. Either way associativity still holds.}.
\end{enumerate}

\section{D-Branes}
\label{subsec:D-branes}
We complete this miscellany of modifications by reintroducing $\Lambda$, the set of \textbf{D-branes}. Objects in the category are now quadruples $(C,O,s,t)$, where $s,t : \{0,1,...,O-1\}\to \Lambda$. An OC Riemann surface between $(C_-,O_-,s_-,t_-)$ and $(C_+,O_+,s_+,t_+)$ is an OC Riemann surface whose free boundaries are labeled by $\Lambda$. The \textbf{free boundaries} are the connected components of the unparameterized boundary\footnote{This should not be confused with the unrelated notion of ``free boundary arcs'' encountered earlier in the context of the reflection principle.}. We require that every open boundary $o \in O_- \cup O_+$ begin at a free boundary with label $s_{\pm}(o)$ and end with at a free boundary labeled $t_{\pm}(o)$.

When composing morphisms, we require that the resultant surface will be labeled in a manner compatible with the two pieces (note that this will always be possible since $s$ and $t$ restrict the composable morphisms).

Morphisms of OC-surfaces (and in particular the QC isomorphisms) are required to respect the labeling of the boundaries. This means that there are more types (i.e., connected components in the moduli space), but the essential construction remains unchanged.

%\section{Enriching By a Local System}
%We now want to turn Segal's category into a category enriched over topological spaces with local systems. For this we first need to define the monoidal category of spaces with local systems, and then enrich the morphism spaces and composition maps accordingly.
%
%\section{Spaces with Local Systems}

\chapter{CFT's and Topological CFT's}
\label{chapt8}
\label{sec:cft-and-tcft}

The category $\mm_\Lambda$ is a symmetric monoidal category. A conformal field theory is a symmetric monoidal functor from this category to the category of vector spaces. These terms will be defined below. We then turn to the definition of topological conformal field theories.

\section{Symmetric Monoidal Categories and Functors}
\label{subsec:monoidal}

We recall the basic definitions of symmetric monoidal categories and functors. For more details we refer the reader to \citep{MacLane}.

\textbf{A monoidal category} consists of the following data:
\begin{itemize}
\item A category $\cl$.
\item A functor $\otimes : \cl \times \cl \to \cl$ with a natural associativity isomorphism $\alpha_{X,Y,Z} : (X \otimes Y) \otimes Z \to X \otimes (Y \otimes Z)$.
\item An identity object $I$ together with natural isomorphisms $\lambda_X : I \otimes X \to X$ and $\rho_X : X \otimes I \to X$.
\end{itemize}
The natural transformations are required to satisfy the ``coherence property'': we want any kind of diagram involving the natural transformations, the tensor operation, and identity morphisms to commute. By a famous theorem of Mac Lane it is sufficient to assume this for only two diagrams. Namely we require that the two ways of transforming $X \otimes (Y \otimes (Z \otimes W))$ to $((X \otimes Y) \otimes Z) \otimes W$ are the same, and that the two ways of transforming $(X \otimes I) \otimes Y$ to $X \otimes Y$ are the same (it may be instructive to write down the corresponding commutative diagrams).

A monoidal category is called \textbf{symmetric} if there exists a natural ``commutativity'' isomorphism $\sigma_{A,B} : A \otimes B \to B \otimes A$, subject to coherence . Again we find that a few relations imply all the rest. Namely, there are now two ways of getting from $A \otimes I$ to $A$, and we require them to be the same. So should the two ways of getting from $(X \otimes Y) \otimes Z$ to $X \otimes (Z \otimes Y)$. Finally, we require that $\sigma_{Y,X} \circ \sigma_{X,Y} = \mbox{Id}_{X,Y}$ (these are the two ways of getting from $X \otimes Y$ to itself).

A \textbf{monoidal functor} between two monoidal categories $\cl, \mathcal{D}$ is a functor $F : \cl \to \mathcal{D}$ together with a natural transformation (not necessarily an isomorphism) $\phi_{X,Y} : F(X) \otimes_{\mathcal{D}} F(Y) \to F(X \otimes_\cl Y)$ and a morphism $\iota : I_{\mathcal{D}} \to F(I_\cl)$, which satisfy the coherence property with respect to all commutative diagrams which may now involve $F$ and $\iota$. The relations are generated by three diagrams: one for $\alpha$, one for $\lambda$, and one for $\rho$. A monoidal functor is \textbf{symmetric} if it satisfies an additional diagram corresponding to $\sigma$ (and then, of course, it is coherent with respect to all such diagrams).

Here are some examples of symmetric monoidal categories which we will run into. We denote only the binary functor and the identity object, the natural isomorphisms being implied. In the future, we will denote only the underlying category, and refer to this list. Here and in what follows, $\kk$ represents a field of characteristic zero.
\begin{enumerate}
\item $(\mbox{Vect}_\kk,\otimes,\kk)$, the category of vector spaces over $\kk$ with the usual tensor product.
\item $(\mbox{Comp}_\kk, \otimes, \kk)$, the category of graded chain complexes over $\kk$. Here $(A \otimes B)_k = \oplus_{r+s = k} A_r \otimes B_s$ and $d(a \otimes b) = da \otimes b + (-1)^p a \otimes db$.
\item $(\mbox{Top}, \times, \mbox{pt})$, topological spaces with the topological product, and a single point.
\item $(\mm_\Lambda, \coprod, \emptyset)$. This is in fact a symmetric monoidal category enriched over topological spaces.
\end{enumerate}

A \textbf{conformal field theory} is a symmetric monoidal functor from $\mm_\Lambda$ to $\mbox{Vect}_\kk$.

\section{Being Rich and Getting Richer}
\textbf{Topological conformal field theories} are symmetric monoidal functors of $C_*(\mm_\Lambda)$, a category obtained from $\mm_\Lambda$ by replacing the morphism spaces with their singular chains. Although we have been informally speaking about an enriched category all along, in order to explain the details of the replacement we now need to be more precise.

Let $\cl$ be a monoidal category. \textbf{A category enriched over} $\cl$, $\mathcal{R}$, is \begin{itemize}
\item A collection of objects $\mbox{Ob}(\mathcal{R})$.
\item For every pair $A,B \in \mbox{Ob}(\mathcal{R})$ an object of $\cl$, $\mathcal{R}(A,B)$.
\item For every $A \in \mbox{Ob}(\mathcal{R})$ a morphism $\epsilon_A : I \to \mathcal{R}(A,A)$.
\item For every triple $A,B,C \in \mbox{Ob}(\mathcal{R})$ a map $\circ_{A,B,C} : \mathcal{R}(A,B) \otimes \mathcal{R}(B,C) \to \mathcal{R}(A,C)$.
\end{itemize}
the usual axioms for a category (associativity, left and right identity) translate into commutative diagrams involving the natural isomorphisms of $\cl$. It is now easy to check that a usual category is a category enriched over sets with cartesian product, and that $\mm_\Lambda$ is enriched over topological spaces with the topological product.

\vspace{2mm}
Let $\cl$ and $\mathcal{D}$ be monoidal categories, and let $(F : \cl \to \mathcal{D}, \phi)$ be a monoidal functor. If $\mathcal{R}$ is a category enriched over $\cl$ we can construct a category $\mathcal{S}$ enriched over $\mathcal{D}$ by the following general procedure
\begin{itemize}
\item $\mbox{Ob}(\mathcal{S}) := \mbox{Ob}(\mathcal{R})$
\item $\mathcal{S}(A,B) := F(\mathcal{R}(A,B))$
\item $(\epsilon_A)_\mathcal{S} := F((\epsilon_A)_\mathcal{R}) \circ \iota$
\item $(\circ_{A,B,C})_\mathcal{S} := F((\circ_{A,B,C})_\mathcal{R}) \circ \phi_{\mathcal{R}(A,B), \mathcal{R}(B,C)}$
\end{itemize}

A straightforward diagram chase shows that $\mathcal{S}$ is indeed a category enriched over $\mathcal{D}$.

\section{From Top to Comp$_\kk$}
We want to replace the morphism spaces by chain complexes that compute their homology. By the previous chapter, this requires a symmetric monoidal functor $\mbox{Top} \to \mbox{Comp}_\kk$. We will show that $C_{\mbox{sing}}$, the singular chains functor over $\kk$, is such a functor if we augment it with the appropriate natural transformations.

For $\phi_{X,Y}$ we take the natural ``shuffle cross product map''
\[(\phi_{X,Y})_n : \oplus_{i = 0}^n C_i(X) \otimes C_{n-i}(Y) \to C_{n}(X\times Y)\]
To construct it, we identify $\Delta^i \times \Delta^j$ with the subset
\[\{(x_0,\hdots,x_i,y_0,\hdots,y_j\} \in \rr^{m+n+2} | 0 \leq x_0 \leq x_1 \leq \cdots \leq x_i \leq 1, \; \; 0 \leq y_0 \leq y_1 \leq \cdots \leq y_j \leq 1\},\]
and divide it into simplices by the various completions of this partial order to a total order of the $x$'s and $y$'s. We map the pair of maps $f : \Delta^i \to X$ and $g : \Delta^j \to Y$ to the singular chain corresponding to $f \times g : \Delta^{i + j} \to X \times Y$, obtained by summing the restrictions to the simplices in the division of $\Delta^i \times \Delta^j$, with the signs determined by the induced orientation from $\rr^{i + j + 2}$. We then extend this definition linearly to arbitrary chains in the tensor product. The details of this construction, as well as analysis of its precise affect on homology (``The Kunneth Formula'') can be found in \citep[pg. 278]{hatcher}.

Note that $C(\mbox{pt})$ is canonically isomorphic to $\kk$, and this defines $\iota$. The coherence conditions for the left and right identity maps and the twist are trivial. To check associativity, it is enough to examine a triplet of maps $f : \Delta^i \to X,\; g : \Delta^j \to Y,\; h : \Delta^k \to Z$. Note that along either path of the diagram

\vspace{5mm}
\xymatrix{
    (CX\otimes CY)\otimes CZ\ar[d]_{\phi_{X,Y}\otimes 1}\ar[r]^{\alpha_D}&CX\otimes(CY\otimes CZ)\ar[d]^{1\otimes\phi_{Y,Z}}\\
    C(X\times Y)\otimes CZ \ar[d]_{\phi_{X \times Y,Z}}&CX \otimes C(Y\times Z)\ar[d]^{\phi_{X,Y\times Z}}\\
    C((X \times Y)\times Z)\ar[r]_{C\alpha_C}&C(X\times(Y\times Z))
}

\noindent
we arrive at the same decomposition of $f \times g \times h$ as a sum of singular simplices, each corresponding to a total ordering of $x_0,\hdots, x_i, y_0,\hdots,y_j,z_0,\hdots,z_k$ extending the partial order $x_0 \leq \cdots \leq x_i,y_0 \leq \cdots$. The orientation of the simplices is induced from the orientation of $\rr^{i + j + k +3}$ in both cases, so the signs are also the same.

\vspace{3mm}
This completes our definition of topological conformal field theories. Two caveats are in order.
\begin{enumerate}
\item To accomodate the most important geometric TCFT's, it is sometimes necessary to add a local system to the morphism spaces and calculate the chains with respect to that local system. This can be thought of as allowing certain projective representations of the category.
\item For technical reasons related to cellular approximation, Costello introduces a variation of the usual singular chains functor which uses regular CW-complexes, rather than simplices, as models.
\end{enumerate}
For the details we refer the reader to \citep{costello}.

\appendix
\chapter{The Type of an Open-Closed-Smooth surface}
\label{AppA}
\label{sec:OCsmoothtype}
Recall that two OCR-surfaces $(C_-,O_-) \Rightarrow (C_+, O_+)$ are of the same OC-smooth type $\tau : (C_-,O_-) \Rightarrow (C_+, O_+)$ if there's a diffeomorphism rel-$\partial$ between them. If we relax the restriction on the boundary, then it is well-known that connected smooth oriented surfaces with boundary are classified up to a diffeomorphism by their genus $g$ and number of boundary components $n$. The following theorem shows that the boundaries are flexible in the smooth category.

\begin{theorem}
Let $\Sigma$ be a connected smooth oriented surface of genus $g$ and with $n$ boundary components. Let $f : \partial \Sigma \to \partial \Sigma$ be an orientation-preserving diffeomorphism of the boundary. Then there's an orientation-preserving diffeomorphism $\overline f : \Sigma \to \Sigma$ with $\overline f_{|\partial \Sigma} = f$.
\end{theorem}

\begin{proof}

The proof hinges on the following technical result:
\begin{lemma}
\label{lemma:smooth twist}
Let $\varphi$ be an orientation-preserving diffeomorphism of $S^1$. Let $A \subset \rr^2$ be an annulus described by two concentric circles. Denote the inner circle of $A$ by $S^1$ and the disc inscribed in it by $D$. Then there's a diffeomorphism $\overline \varphi : \rr^2 - D$ which is the identity on the complement of $A$, and such that $\overline \varphi_{|S^1} = \varphi$. Furthermore, if $\varphi$ is rotation by a constant angle we can extend $\overline \varphi$ to a diffeomorphism of $\rr^2$ which acts by rotation on $D$.
\end{lemma}

\begin{proof}
To show this, it suffices to find a smooth path $\varphi_t$ through the group of diffeomorphisms of $S^1$ such that $\varphi_0 = \varphi$ and $\varphi_1 = id$, and with
\[\frac{\partial^n}{\partial t^n}_{|t = 0} \varphi_t(x)= \frac{\partial^n }{\partial t^n}_{|t = 1}\varphi_t(x) = 0\] for all $x \in S^1$ and $n \geq 1$.

Indeed, using such a path we can define $\overline \varphi$ on $A$ in the obvious way: $\overline \varphi(r,\theta) := (r,\varphi_{t(r)} (\theta))$, where $t(r) = \frac{r - r_1}{r_2 - r_1}$. One easily verifies that the Jacobian of this map does not vanish and that the extension by identity to $\rr^2 - D$ is smooth at the outer boundary circle of $A$. In fact we can always extend $\overline \varphi$ to $\rr^2 - \{0\}$ by letting $\overline \varphi(r,\theta) := (r, \varphi(\theta))$ in $D$. This will be smooth at $0$ if $\varphi$ is rotation by a constant angle.

We now construct such a path. Let $\sigma(t) : [0,1] \to [0,1]$ be a smooth surjective map with all derivatives vanishing at $0$ and $1$. Assume first that $\varphi$ fixes the basepoint. Then we can think of $\varphi$ as a smooth map $[0,2\pi] \to [0,2\pi]$ with $\varphi'(x) > 0$, $\varphi(0) = 0$ and $\varphi(2\pi) = 2\pi$. For every $0 \leq t \leq 1$ the map $\varphi_t(x) = (1-\sigma(t)) \varphi(x) + \sigma(t)x$ enjoys the same properties and is therefore a diffeomorphism of $S^1$. This clearly gives a smooth path as desired.

If $\varphi$ does not fix the basepoint $x_0$ we define the path on $[0,1/2]$ by $e^{i \delta \sigma(2t)}\varphi(x)$, where $\delta$ is the angle between $x_0$ and $\varphi(x_0)$. We can then smoothly extend this path to $[1/2,1]$ in the same way as before. Note that the path is smooth also at $t = 1/2$ because of the vanishing of the derivatives in $t$.
\end{proof}

We return to the proof of the theorem.

By the classification theorem of smooth surfaces we may assume $\Sigma$ is the standard sphere with $g$ handles attached to its southern hemisphere and $n$ standard discs removed from its northern hemisphere at prescribed locations. For example, we may identify the open hemisphere with $\rr^2$ and assume the discs are of radius 1/10 and centered at $(k,0)$ for $k = 1,...,n$. The precise details are immaterial, the point is that we can write very explicit formulae for diffeomorphisms which are the identity outside of some compact set, and which realize any desired permutation of the discs. For example, we can swap two neighbouring discs by a half-rotation of a disc $D$ containing them. Applying the lemma, we may extend this diffeomorphism of $D$ to all of $\rr^2$ so that it is the identity outside some disc containing the two swapped discs. If we take this disc small enough it will not intersect any other discs.

After we obtain the desired permutation of the boundary circles, we may appeal to the lemma once again to obtain the desired map on each boundary circle without affecting the diffeomorphism near other boundary components: we take an annulus whose inner circle is the boundary circle and which does not intersect any other boundary circles.

\begin{figure}[hb]
  \centering
\includegraphics[width = 40mm]{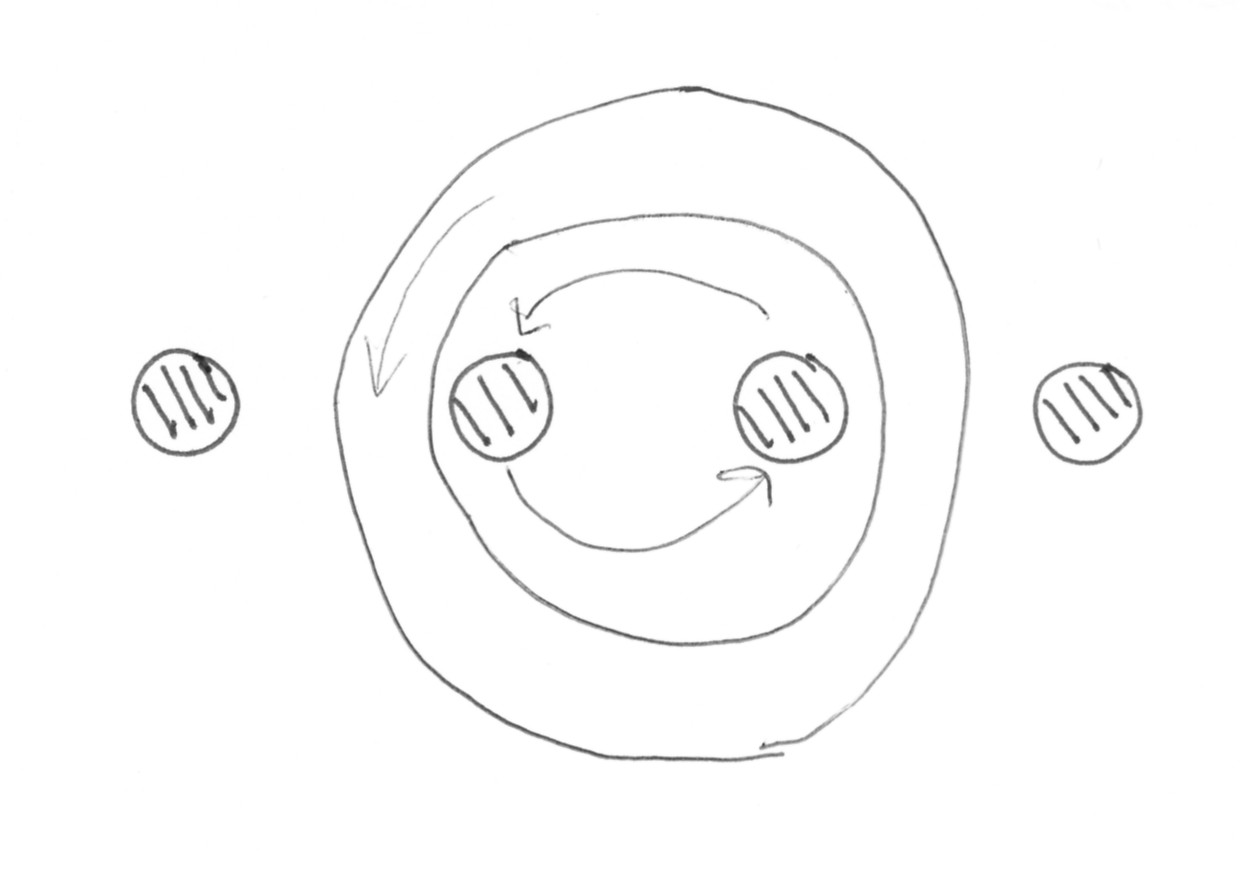}\hspace{10mm}
\includegraphics[width = 40mm]{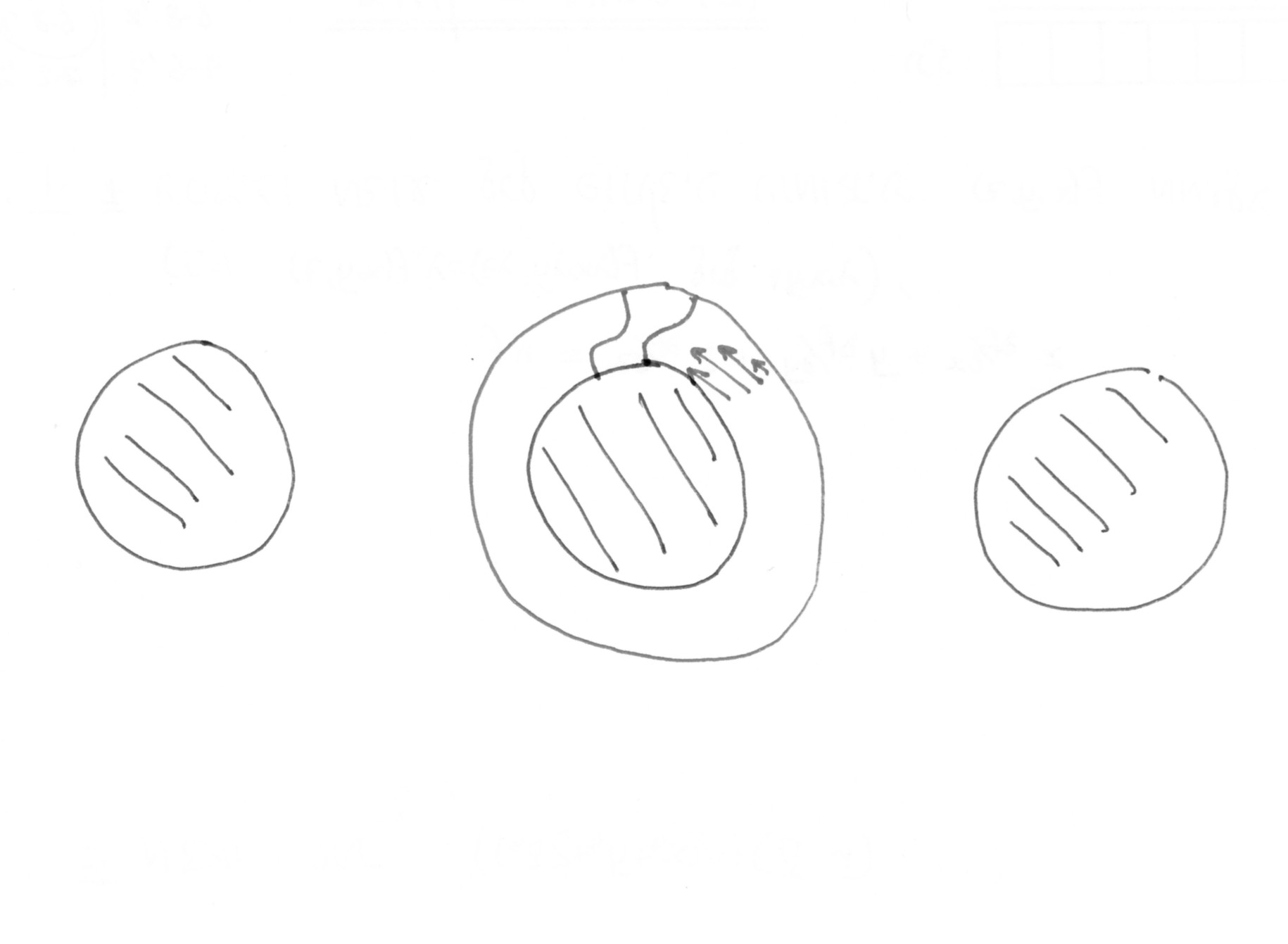}
  \caption[aligning]%
  {Using diffeomorphisms we can swap boundary components (left), and also obtain any desired smooth reparametrization (right).}
\end{figure}

\end{proof}

\begin{corollary}
\label{cor:type-data}
The OC-smooth type of a surface is determined uniquely by the following (discrete) data:
\begin{itemize}
\item A partition of $\{0,1,...,O_+ -1\} \cup \{0',1',...,(O_- -1)'\}$ into cycles $\omega^1,\omega^2,...,\omega^r$ describing the order in which the open boundary intervals are traversed when moving along boundary circles with the outward normal to your left.
\item An unordered list of ``component data'': for each path component, what is its genus, how many boundary components does it have, and which closed boundary components and cycles of open boundaries map to it.
\end{itemize}
\end{corollary}

\begin{proof}[Proof Sketch.]
Note that if  $B, B' : [0,1]^{n_+} \times [1,0]^{n_-} \to S^1$ are smooth orientation-preserving embeddings ($[1,0]$ represents the unit interval with reverse orientation) then there is an orientation-preserving diffeomorphism $\gamma : S^1 \to S^1$ extending $B' \circ B^{-1}$ iff the cyclic order of the embedded intervals is the same.

The corollary now follows from the classification theorem together with Theorem 11.1 above.
\end{proof}

\begin{remark}
\label{rem:composing-types}

Although an OC-smooth structure ``contains more information'', in the sense that there are many OC-smooth structures compatible with any given OC-QC structure, the isomorphism \emph{type} of OC-quasiconformal and OC-smooth surface are the same - in fact they both reduce to the topological type (cf. \textsection \ref{subsec:the morphism space}). From the universal property of lemma \ref{lemma:pushout} we see that there is a unique type associated to the pushout of surfaces of types $\tau_1, \tau_2$, denoted by $\tau_2 \circ \tau_1$. Indeed, it is not hard to derive an explicit recipe for its calculation:
\begin{itemize}
\item Path components of $\tau$ are equivalence classes of path components of $\tau_1$ and $\tau_2$, where the equivalence relation is generated by identifying components which share a closed or an open boundary component.
\item The euler characteristic of each path component is calculated by $\chi(A \cup B) = \chi(A) + \chi(B) - \chi(A \cap B)$, where $\chi(A \cap B)$ is the number of common open boundary components.
\item Open boundary cycles are spliced and reglued according to the common open boundary intervals. This determines the number of boundary circles of each component, the partition of $\{0,...,O_--1,0',...,(O_+-1)'\}$ into cycles, and which cycle maps into which component.
\item The genus of each component is calculated from the euler characteristic and number of boundary components.
\end{itemize}
\end{remark}

\chapter{Smooth Sewing of Complex Structures}
\label{AppB}
\label{sec:smooth sewing}
We show that if we work with smooth parameterizations, the sewed complex structure agrees with the smooth structures of the pieces up to the boundary.

More precisely, let $\overline \hh$ and $\overline \hh^*$ be the closed upper and lower half planes.

\begin{theorem}
Given a function $\rho : \rr \to \rr$ whose extension $\tilde \rho : S^1 \to S^1$ is smooth there exists a (unique) complex structure on $\cc = \overline \hh \cup_\rho \overline \hh^*$ such that a function $f : U \to \cc$ is analytic iff $f_{| \overline \hh}$ and $f_{|\overline \hh^*}$ are smooth and, away from $\rr$, analytic. In particular, the image of the real axis under any conformal chart is smooth.
\end{theorem}

\begin{proof}
Denote by $j$ the standard complex structure on $\cc$.

We will construct a sequence of diffeomorphisms $\tau_k : U_k \to V_k$, where $U_k,V_k$ are open neighbourhoods of $\rr$ in $\overline \hh$, such that
\begin{enumerate}[(a)]
\item $\tau_{k|\rr} = \rho$
\item $\tau_k$ is QC (i.e., its dilatation is bounded), and
\item $J_k := \tau_{k*} (j_{|\overline \hh})$ and $j_{|\overline \hh^*}$ agree to order $C^k$ on $\rr$.
\end{enumerate}

This will involve repeated iteration of the exponential map. Before we go into the details of this construction, let us explain why this gives the desired result. It is well known\footnote{see \citep{bers}, pg. 1084. The requirement that $\mu$ be defined in the entire complex plane is clearly inessential as we can apply a smooth extension operator.} that given an almost complex structure $J$ of bounded dilatation on an open $W \subset \cc$ which is $C^k$, there exists a conformal mapping $\omega : (W,J) \to (W',j)$, $W' \subset \cc$, \emph{and this map is $C^k$}\footnote{cf. theorem \ref{thm:integration}}. Furthermore, this map is unique up to conformal mappings of $W'$. Taking $J = j_{|\overline \hh^*} \cup J_k$, defined on $W = \overline \hh^* \cup U_k$, we see that for arbitrary $k$'s we can produce a $C^k$ chart for the sewed complex structure in a neighbourhood of every $z \in \cc$. By uniqueness, all of these charts are conformally equivalent and so we see that, in-fact, all these charts define the same $C^\infty$ structure on $\cc$ which agrees with the smooth structures on the closed half planes. The theorem follows readily from this.

We now construct the diffeomorphisms $\tau_k$. We take $U_{-1} = \overline \hh$, and define $\tau_{-1}$ by $\tau_{-1}(x,y) = (\rho(x), y)$. Let
$\hat x \equiv \left(
                 \begin{array}{c}
                   1 \\
                   0 \\
                 \end{array}
               \right)$
denote the constant unit vector field in the $x$-direction. We have $J_{-1} \hat x = (d\tau_{-1} \circ j \circ \left(d\tau_{-1}\right)^{-1}) \hat x = \left(
                 \begin{array}{c}
                   0 \\
                   \rho'(x) \\
                 \end{array}
               \right)$.

By elementary calculus, we can define a diffeomorphism $\delta_0$ in an open neighbourhood $U_0 \subset \overline \hh$ of $\rr$ by
\[\delta_0(\gamma_x(y)) = (x,y)\]
where $\gamma_x$ is the unique integral curve of the vector field $J_{-1} \hat x$ which satisfies $\gamma_x(0) = x$. We have $\delta_{0|\rr} = id$, $(\delta_0)_* \left(J_{-1} \hat x\right) = \hat y$ (everywhere in $U_0$) and on the real axis we have $(\delta_0)_* \hat x = \hat x$. Setting $\tau_0 = \delta_0 \circ \tau_{-1}$ we see that $\tau_0$ satisfies conditions (a) - (c) above.

We define the following $\tau_k$ recursively in the same fashion: let $\delta_k$ be the diffeomorphism defined in an open neighbourhood of $\rr$, $U_k$, by the equation $\delta_k(\gamma_x(y)) = (x,y)$, for $\gamma_x$ the integral solution of $J_{k-1} \hat x$ which satisfies $\gamma_x(0) = x$. $\tau_k = \delta_k \circ \tau_{k-1}$.

We prove that $\tau_k$ satisfies condition (c) by induction on $k$. A smooth, real-valued function $f$ defined in the neighbourhood of the real axis is said to belong to $\Theta(y^k)$ if $\lim_{y \to 0} \frac{f(x,y)}{y^k}$
exists for every $x \in \rr$. We will write $\Theta(y^k)$ in place of such a function, and apply the obvious rules for adding or multiplying such classes.

Assume that $J_{k-1} \hat x = \left(
                 \begin{array}{c}
                   \Theta(y^m) \\
                   1 + \Theta(y^n) \\
                 \end{array}
               \right)$
(for $J_{-1}$, this holds with $m = \infty$ and $n = 0$). By the mean value theorem we have
\[\gamma_u(v) = \left(
                 \begin{array}{c}
                   u + v\Theta(v^{m}) \\
                   v + v\Theta(v^{n}) \\
                 \end{array}
               \right) = \left(
                 \begin{array}{c}
                   x  \\
                   y  \\
                 \end{array}
               \right)\]
after taking the inverse mapping we obtain
\[\delta_k\left(
                 \begin{array}{c}
                   x  \\
                   y  \\
                 \end{array}
               \right) = \left(
                 \begin{array}{c}
                   u \\
                   v \\
                 \end{array}
               \right) = \left(
                 \begin{array}{c}
                   x + \Theta(y^{m+1}) \\
                   y + \Theta(y^{n+1}) \\
                 \end{array}
               \right)\]
and
\[D\delta_k = \left(
                \begin{array}{cc}
                  1 + \Theta(y^{m+1}) & \Theta(y^{m}) \\
                  \Theta(y^{n+1}) & 1 + \Theta(y^{n}) \\
                \end{array}
              \right)\]

$J_k$ is the image of $J_{k-1}$ under $\delta_k$, and so $\delta_{k*}\hat x \overset{J_k} \mapsto  \delta_{k*}J_{k-1} \hat x = \hat y$. The unique almost complex structure that maps $\left(
                 \begin{array}{c}
                   A \\
                   B \\
                 \end{array}
               \right)$ to $\left(
                 \begin{array}{c}
                   0 \\
                   1 \\
                 \end{array}
               \right)$ is given by $\left(
                                      \begin{array}{cc}
                                        B & A \\
                                        \frac{1 + B^2}{A} & B \\
                                      \end{array}
                                    \right)$ and so we find that
\[J_k \hat x = \left(
                 \begin{array}{c}
                   \Theta(y^{n+1}) \\
                   \frac{1 + \Theta(y^{2n + 2})}{1 + \Theta(y^{m+1})} \\
                 \end{array}
               \right) = \left(
                 \begin{array}{c}
                   \Theta(y^{n+1}) \\
                   1 + \Theta(y^{n'}) \\
                 \end{array}
               \right)\] where $n' = \min\left(2n + 2,m+1\right)$.

It is now a simple matter to check that starting with $m = \infty, n =0$ we have $m,n \to \infty$ (proof by induction; except at the first step one finds that $n' = m +1$).

This completes the proof that $\tau_k$ satisfy property (c). By restricting to a smaller open neighbourhood whose closure is contained in $U_k$, we may ensure that the dilatation is bounded. Property (a) obviously holds as well. This completes the proof of the theorem.
\end{proof}

\newpage
\bibliographystyle{plainnat}
\bibliography{thesis_bibliography}
\nocite{*}

\end{document}